\documentclass[11pt]{amsart}

\usepackage{amssymb,amsthm,amsmath}
\usepackage[numbers,sort&compress]{natbib}
\usepackage{color}
\usepackage{graphicx}
\usepackage{tikz}
\usepackage[colorlinks,
linkcolor=red,
anchorcolor=blue,
citecolor=green
]{hyperref}

\def\a{\alpha}

\def\d{\delta}

\def\e{\varepsilon}

\let\newpf\proof \let\proof\relax 
\newenvironment{pf}{\newpf[\proofname]}{\qed\endtrivlist}

\newcommand{\ba}{\overline{A}}

\def\be{\begin{equation}}
\def\ee{\end{equation}}

\def\ba{{\begin{align}}}
\def\ea{{\end{align}}}

\def\bm{\begin{matrix}}
\def\em{\end{matrix}}

\def\a{{\alpha}}

\def\d{{\underline d}}

\def\0{{\mathbf 0}}

\newtheorem{Theorem}{Theorem}[section]
\newtheorem{Lemma}{Lemma}[section]
\newtheorem{Proposition}{Proposition}[section]
\newtheorem{Corollary}{Corollary}[section]
\newtheorem{Remark}{Remark}[section]

\newtheorem{Definition}{Definition}[section]

\numberwithin{equation}{section}

\theoremstyle{definition}

\renewcommand{\mod}{\operatorname{mod}}

\newcommand{\C}{{\mathbb C}}

\newcommand{\Q}{{\mathbb Q}}
\newcommand{\R}{{\mathbb R}}
\newcommand{\T}{{\mathbb T}}

\newcommand{\Z}{{\mathbb Z}}

\def\B0{{\bold{0}}}

\catcode`\@=12

\def\Empty{}
\newcommand\oplabel[1]{
  \def\OpArg{#1} \ifx \OpArg\Empty {} \else
    \label{#1}
  \fi}

\newcommand{\comm}[1]{}
\newcommand{\comment}[1]{}

\begin{document}

\title[]{Arithmetic version of Anderson localization via reducibility}

\author{Lingrui Ge}
\address{
Department of Mathematics, University of Califoria Irvine, CA, 92697-3875, USA}

 \email{lingruige10@163.com}

\author {Jiangong You}
\address{
Chern Institute of Mathematics and LPMC, Nankai University, Tianjin 300071, China} \email{jyou@nankai.edu.cn}

\begin{abstract}
The arithmetic version of Anderson localization (AL), i.e.,  AL with explicit arithmetic description on both  the localization frequency and   the localization phase, was first given by Jitomirskaya \cite{J} for the almost Mathieu operators (AMO). Later, the result was generalized by Bourgain and Jitomirskaya  \cite{bj02} to a class of {\it one dimensional} quasi-periodic long-range operators.
In this paper, we propose a novel approach based on  an arithmetic version of Aubry duality and quantitative reducibility.
Our method enables us to prove the same result for the class of quasi-periodic long-range operators in {\it all dimensions}, which includes \cite{J, bj02} as special cases.
\end{abstract}

\maketitle

\section{Introduction}
In this paper, we consider the quasi-periodic long-range operator on $\ell^2(\Z^d)$:
\begin{equation}\label{1.2}
(L_{V,\alpha,\theta}u)_n=\sum\limits_{k\in\Z^d}\hat{V}_ku_{n-k}+2\cos 2\pi(\theta+\langle n,\alpha\rangle)u_{n}, \ \ n\in \Z^d,
\end{equation}
where $V(x)=\sum_{k\in\Z^d}\hat{V}(k)e^{2\pi i\langle k,x\rangle}\in C^r(\T^d,\R)$ $(r=0,1,\cdots,\infty,\omega)$, $\theta\in \T$ is called the phase and $\alpha\in\T^d$ is called the frequency.

Operator \eqref{1.2} has received a lot of attentions \cite{aj1,bj02,gyz,bj,cd,jk2,b1} since the 1980s. On  one hand, the spectral properties of operator \eqref{1.2} have close relation to that of its Aubry dual
\begin{equation}\label{1.1}
(H_{V,\alpha,x}u)_n=u_{n+1}+u_{n-1}+V(x+n\alpha)u_{n},\ \  n\in \Z.
\end{equation}
For partial references, one may consult \cite{bj02,aj1,ayz,ayz1,gyz,aj,jk2,p,gjls}. On the other hand, operator \eqref{1.2} itself contains several popular quasi-periodic models. If we take $V(x)=\sum_{i=1}^d2\lambda^{-1}\cos2\pi{x_i}$,  \eqref{1.2} is reduced to  quasi-periodic Schr\"odinger operator on $\ell^2(\Z^d)$:
\begin{equation}\label{highschrodinger}
H_{\lambda,\alpha,\theta}= \Delta +  2\lambda\cos2\pi(\theta+ \langle n,\alpha\rangle)\delta_{nn'},
\end{equation}
where $\Delta$ is the usual Laplacian on $\Z^d$ lattice. If $d=1$, operator \eqref{highschrodinger} is the famous almost Mathieu operator (AMO).

{\it Anderson localization} (pure point spectrum with exponentially decaying eigenfunctions) for quasi-periodic operators has been widely studied over the past sixty years. We give a brief review here. As one can see, both the operators \eqref{1.2} and \eqref{1.1} are  a family of operators parameterized by two parameters,  $(\alpha,\theta)$ in \eqref{1.2} and  $(\alpha,x)$ in \eqref{1.1}, when $V$ is fixed. When $d=1$, Fr\"ohlich-Spencer-Wittwer \cite{fsw} and Sinai \cite{Sinai} proved that $H_{\lambda V,\alpha,x}$ has Anderson localization  for $a.e.$ $x$ and large enough $\lambda$ if $V$ is cosine-like and $\alpha$ is  Diophantine \footnote{ $\alpha \in\T^d$ is  called {\it Diophantine}, denoted by $\alpha \in {\rm DC}_d(\kappa,\tau)$, if there exist $\kappa>0$ and $\tau>d-1$ such that
\begin{equation}\label{dio}
{\rm DC}_d(\kappa,\tau):=\left\{\alpha \in\T^d:  \inf_{j \in \Z}\left| \langle n,\alpha  \rangle - j \right|
> \frac{\kappa}{|n|^{\tau}},\quad \forall \  n\in\Z^d\backslash\{0\} \right\}.
\end{equation}
Let ${\rm DC}_d:=\bigcup_{\kappa>0} {\rm DC}_d(\kappa,\tau)$.}. Eliasson \cite{Eli97} proved that if $V$ is a Gevrey function satisfying non-degenerate conditions, for any {\it fixed} Diophantine $\alpha$, $H_{\lambda V,\alpha,x}$ has pure point spectrum for $a.e.$ $x$ and large enough $\lambda$.   Bourgain and Goldstein \cite{bg} proved that, in the positive Lyapunov exponent regime, for any {\it fixed} $x$, $H_{\lambda W,\alpha,x}$ has AL for  $a.e.$  Diophantine $\alpha$ provided that  $V$ is a non-constant real analytic function. Bourgain and Jitomirskaya \cite{bj} generalized the result in \cite{bg} to certain band models. Klein \cite{K} generalized the results in \cite{bg} to more general Gevrey potentials.

We would like to remark that \cite{fsw,Sinai,Eli97} gave an arithmetic description on the localization frequency $\alpha$, but they didn't have an arithmetic description on the localization phase $\theta$, in contrast, \cite{bg,bj,K} proved localization  for fixed phase $\theta$, but there is no arithmetic description on the localization frequency $\alpha$. Thus, for any fixed $(\alpha,x)$, it is not clear whether $H_{\lambda W,\alpha,x}$ has AL. The breakthrough belongs to Jitomirskaya \cite{J} who proved  that for any {\it fixed} Diophantine frequency and {\it any  fixed} Diophantine phase\footnote{ $\theta \in\T$ is  called {\it Diophantine}, denoted by $\alpha \in {\rm DC}_\alpha(\kappa,\tau)$, if there exist $\kappa>0$ and $\tau>d-1$ such that
\begin{equation}
{\rm DC}_\alpha(\kappa,\tau):=\left\{\theta \in\T:  \inf_{j \in \Z}\left| \langle n,\alpha  \rangle - 2\theta-j \right|
> \frac{\kappa}{(|n|+1)^{\tau}},\quad \forall \  n\in\Z^d \right\}.
\end{equation}
Let $\Theta:=\bigcup_{\kappa>0} {\rm DC}_\alpha(\kappa,\tau)$.}, the almost Mathieu operator has AL if $|\lambda|>1$. Recently, such arithmetic description on the frequency and the phase was explored even in a sharp way \cite{JLiu,JLiu1}. Namely, for Diophantine phase, there is a sharp spectral transition in frequency \cite{JLiu}. For Diophantine frequency, there is a sharp spectral transition in phase \cite{JLiu1}.
For more general quasi-periodic long-range operators \eqref{1.2}, Bourgain and Jitomirskaya \cite{bj02} generalized the methods in \cite{J} and proved that for any {\it fixed} Diophantine frequency and {\it fixed} Diophantine phase, $L_{\lambda V,\alpha,\theta}$ has AL for small $\lambda$ and real analytic $V$. For more such kinds of results, one may consult \cite{aj1} e.t.c. The localization results which have  arithmetic descriptions on both the frequency and the phase, are called {\it arithmetic version of Anderson localization}.

We emphasize that all the above localization results are restricted to one dimensional quasi-periodic operators. If the dimension $d>1$, one meets essential difficulties and   AL results are quite few. When $d=2$, Bourgain, Goldstein and Schlag \cite{bgs2} proved that  for any {\it fixed} $\theta$, $H_{\lambda,\alpha,\theta}$ has AL for sufficiently large $\lambda$  and a positive Lebesgue measure set of $\alpha$.   As a matter of  fact, we point out that Bourgain-Godstein-Schlag's result \cite{bgs2} is for (\ref{highschrodinger}) with more general potentials
$
V(n_1,n_2)=v(\theta_1+n_1\alpha_1, \theta_2+n_2\alpha_2)
$,
where $v$ is a real analytic function on $\T^2$  which is non-constant on any horizontal or vertical line. Later, Bourgain \cite{b} generalized the result \cite{bgs2} to the case when $d\geq 3$.  Recently,
the results in \cite{b} have been largely extended by Jitomirskaya-Liu-Shi
\cite{jls} to general analytic $k$-frequency quasi-periodic operators on $\Z^d$
for arbitrary $k, d$. Inspired by  pioneer work  of  Fr\"ohlich-Spencer-Wittwer \cite{fsw} and Sinai \cite{Sinai}, Chulaevsky-Dinaburg \cite{cd} proved  essentially the same result as that in \cite{b} for higher dimensional quasi-periodic long-range operator
\[
(Lu)_n=\sum\limits_{k\in\Z^d}\hat{V}_ku_{n-k}+\lambda W(\theta+\langle n,\alpha\rangle)u_{n}, n\in \Z^d,
\]
 with cosine-like potentials $W$,  i.e.,  for any fixed phase $\theta$, the operator  has AL for positive measure $\alpha$ if the coupling constant $\lambda$ is sufficiently large. Finally,
 Bourgain \cite{b1} obtained the same result for the operator $L$ if  $V$ and $W$ are real analytic and $\lambda$ is sufficiently large.
We emphasize that  all the above higher dimensional AL results are for fixed phase and a {\it positive measure} set of frequency (depending on the phase).
These results are weaker than those in the one dimensional case where one can prove AL for fixed phase and a {\it full measure} set of frequency (again, not uniform for all $\theta$).

Recently, a new approach for proving Anderson localization based on Aubry duality and reducibility was given by  Avila,You and Zhou \cite{ayz}, Jitomirskaya and Kachkovskiy \cite{jk2} which leads to stronger conclusions.  Firstly, this method doesn't care about the dimension and the regularity of the potential since reducibility holds for multi-frequency cocycles \cite{ds,e,LYZZ,gyz} and low regularity cocycles \cite{CG,ccyz}. Based on this, Jitomirskaya and Kachkovskiy  \cite{jk2} proved that  for any {\it fixed} Diophantine $\alpha$ and real analytic $V$, $L_{\lambda V,\alpha,\theta}$ has pure point spectrum for small enough $\lambda$ and $a.e.$ $\theta$. Recently, Ge-You-Zhou \cite{gyz} established a general criterion and proved exponential dynamical localization in expectation for $L_{\lambda V,\alpha,\theta}$ under the same condition. Secondly, in one  dimensional case, one can even get non-perturbative localization in a sharp way by non-perturbative reducibility method. This ultimately leads  Avila-You-Zhou's \cite{ayz} to the solution of the measure version of Jitomirskaya's conjecture \cite{j}.

However, if $d>1$,  localization is not known for any concrete  $(\alpha,\theta)$. In other words, there is no arithmetic version of Anderson localization. People even do not know if one should expect such result since higher dimensional case is more ``random" than the one dimensional case. From a methodological point of view, the method in \cite{J,bj02} seems non-trivial to be generalized to the higher dimensional case due to obstacles in computing the Green function. While the methods developed by \cite{ayz,jk2} lose control of a zero Lebesgue measure set of phase which leads to loss of an arithmetic description of the localization phase.

In this paper, by introducing an auxiliary measure defined by reducibility, we find a strategy to recover the phases lost in using the method in \cite{ayz,jk2}. We thus develop an arithmetic-theoretic Aubry duality which gives, for any fixed Diophantine frequency $\alpha$, the one-to-one correspondence between the localization phases of \eqref{1.2} and the reducibility energies (described by the rotation number) of \eqref{1.1}.  Since reducibility theory can provide a clear arithmetic description on the rotation number of the eigenvalue equation of \eqref{1.1}, one thus has a clear arithmetic description on the localization phase of the dual model \eqref{1.2} and proves Anderson localization for operators \eqref{1.2} for all Diophantine frequencies and all Diophantine phases.

 Now we state our results. Assume $V\in C^\omega(\T^d,\R)$, we denote by
$$
S_E^{V}(x):=
\begin{pmatrix}
E-V(x) & -1\\
1 & 0
\end{pmatrix},   \quad E\in\R.
$$
Recall that the cocycle $(\alpha,S_E^V)$ is said to be $C^\omega$-almost reducible if there exists a sequence $B_n\in C^\omega(\T^d,PSL(2,\R))$ such that $B^{-1}_n(x+\alpha)S_E^V(x)B_n(x)$ converges to a constant matrix. Let
\begin{align*}
\mathcal{AR}=\{E\in \Sigma_{\alpha,V}| \mbox{$(\alpha,S_E^V)$ is $C^\omega$-almost reducible}\},
\end{align*}
where $\Sigma_{\alpha,V}$\footnote{See Section 2.3 for more details.} is the spectral set of $H_{V,\alpha,x}$.
\begin{Theorem}\label{long}
Assume $\alpha\in DC_d$, $V\in C^\omega(\T^d,\R)$ and $\Sigma_{\alpha,V}=\mathcal{AR}$. Then $L_{V,\alpha,\theta}$ has Anderson localization for $\theta\in \Theta$, where $\Theta$ is the set defined in footnote 2.
\end{Theorem}

We remark that Theorem \ref{long} is a global result, namely, we don't need to assume that $V$ is small.
However, by Theorem A in \cite{e}, $\mathcal{AR}=\Sigma_{\alpha,V}$ always holds if  $V\in C^\omega(\T^d,\R)$  is sufficiently small. Thus we immediately have the following corollary.
\begin{Corollary}\label{high1}
Assume that $\alpha\in DC_d$. Then there exists $\lambda_0(\alpha,V,d)$, such that
$L_{\lambda V,\alpha,\theta}$ has Anderson localization for $\theta\in\Theta$ if $\lambda<\lambda_0$.
\end{Corollary}
\begin{Remark}
Corollary \ref{high1} is the first higher dimensional Anderson localization result with explicit arithmetic description on both the frequency and the phase. Anderson localization with arithmetic description on the frequency but not on the phase was proved in \cite{gyz,jk2}.
\end{Remark}
\begin{Remark}
If $d=1$, Bourgain and Jitomirskaya \cite{bj02} proved that,  under exactly the same assumption, Anderson localization holds for $\theta\in \Theta$\footnote{The Anderson localization phase given in \cite{bj02} is a little larger than $\Theta$, we mention that we can prove Anderson localization for the same phase set with minor modifications of the proof.}. Thus our result can be regarded as a generalization of their result to higher dimensions.
\end{Remark}

If the dimension is one, by the global theory developed by Avila \cite{avila0}, the spectrum set $\Sigma_{\alpha,V}$ can be decomposed as three regimes, the subcritical regime $\Sigma_{\alpha,V}^{sub}$, the critical regime $\Sigma_{\alpha,V}^{cri}$ and the supcritical regime $\Sigma_{\alpha,V}^{sup}$. Moreover, it was proved by Avila in \cite{avila1,avila2}, $\Sigma^{sub}_{\alpha,V}=\mathcal{AR}$. Thus we can re-state  Theorem \ref{long} as the following corollary.
\begin{Corollary}\label{sub1}
Assume that $\alpha\in DC_1$, $V\in C^\omega(\T,\R)$ and $\Sigma_{\alpha,V}=\Sigma_{\alpha,V}^{sub}$. Then $L_{V,\alpha,\theta}$ has Anderson localization for $\theta\in \Theta$.
\end{Corollary}
\begin{Remark}
Corollary \ref{sub1} is a generalization of Jitomirskaya's localization result in \cite{J} for  supercritical AMO.
\end{Remark}



A key part of our proof is Aubry duality, which has a long history starting from \cite{aubryandre}. As pointed out in \cite{aubryandre,jk2}, Aubry duality can be understood as the correspondence between absolutely continuous (ac) and pure point (pp) spectra of $\{H_{V,\alpha,x}\}_{x\in T^d}$ and $\{L_{V,\alpha,\theta}\}_{\theta\in\T}$ and vice versa. However, such correspondence may not hold, even for the almost Mathieu operator, one can see \cite{L,jk2} for more details.

A more refined notion than absolutely continuous spectrum is the {\it reducibility} of the corresponding Schr\"odinger cocycle $(\alpha, S_E^V)$, see \cite{e}. On the one hand, it is explored in \cite{jk2,ayz} that reducibility of $\{(\alpha, S_E^V)\}_{E\in\R}$ for $a.e.$ $E$ with respect to the integrated density of states $\mathcal{N}$\footnote{See Section 2.2 for definition.} implies pure point spectrum of $\{L_{V,\alpha,\theta}\}_{\theta\in\T}$ for $a.e.$ $\theta$. On the other hand, as proved in \cite{p,aj1}, pure point spectrum of $\{L_{V,\alpha,\theta}\}_{\theta\in\T}$ for $a.e.$ $\theta$ implies reducibility of $\{(\alpha, S_E^V)\}_{E\in\R}$ for $a.e.$ $E$ with respect to the integrated density of states.
The above correspondence can be viewed as a {\it measure-theoretic level of Aubry duality}.

Generally speaking, Aubry duality has been explored at different levels, i.e., the  operator-theoretic level \cite{gjls}, the quantitative level \cite{aj1} and the measure-theoretic level \cite{jk2,ayz}.
Each level of duality has its own important applications. Based on the operator-theoretic duality developed in \cite{gjls}, Jitomirskaya \cite{J} completely solved the measure version of Aubry-Andre conjecture \cite{aubryandre} by non-perturbative localization methods. Based on the quantitative duality, Avila and Jitomirskaya \cite{aj1} proved a non-perturbative Eliasson's almost reducibility theory, gave a sharp H\"older exponent of integrated density of states and solved the ``Dry Ten Martini Problem" for non-critical AMO with Diophantine frequency. Based on the measure-theoretic duality, Avila, You and Zhou \cite{ayz} solved  the measure version of Jitomirskaya's conjecture \cite{j} on sharp phase transitions of AMO by non-perturbative reducibility methods. These applications of different levels of duality imply the general philosophy that Aubry duality is a bridge connecting operators and their duals. Once one knows everything about the operators themselves at
 some level, one knows everything about their duals at this level.

Thus,  an interesting question is whether  one can develop a new level of Aubry duality, which gives the arithmetic version of Anderson localization? This level of duality should completely transfer the arithmetic description of reducibility energy of $\{(\alpha, S_E^V)\}_{E\in\R}$ to the arithmetic description of the localization phase of $\{L_{V,\alpha,\theta}\}_{\theta\in\T}$. More precisely, we hope to give a full measure set $\Theta$ with precise arithmetic description, such that  $\{L_{V,\alpha,\theta}\}_{\theta\in\T}$ has AL for $\theta\in \Theta$ if all $\{(\alpha, S_E^V)\}_{E\in\R}$ with  the rotation number in  $\Theta$ are reducible. We call such kind of duality {\it arithmetic-theoretic Aubry duality}. Before further explanations, we give another reason why arithmetic-theoretic Aubry duality is important. Recently, Jitomirskaya and Liu \cite{JLiu1} further developed their method in \cite{JLiu}, and proved sharp spectral transition in phase for Diophantine AMO. More importantly, they uncover a new type of hierarchy for quasi-periodic operators which is called reflective-hierarchy structure. As pointed out in \cite{JLiu1}, this progress also uncovers some general phenomena for quasi-periodic operators that spectral transition happens not only in frequency, but also in phase. While the spectral transition in phase need a complete arithmetic description on phase which is out of reach by measure-theoretic Aubry duality. This makes an arithmetic version of Aubry duality of particular importance.

We remark that measure-theoretic Aubry duality is much easier than arithmetic-theoretic Aubry duality. The shortcoming is the measure-theoretic Aubry duality gives AL  for almost every phase, and loses control of a zero measure subset in $\Theta$. To prove AL for all $\theta\in \Theta$,  i.e., the arithmetic version of AL, one only need to prove $d\mu^{pp}_\theta$ is continuous in $\Theta$. However, this is a difficult job and we don't know how to prove it since $d\mu^{pp}_\theta$ sensitively depends on $\theta$.   Our strategy is to  define a new measure $d\nu_\theta$ via reducibility, we call it $\mathcal{R}$-measure, which is absolutely continuous with respect to $d\mu^{pp}_\theta$. The advantage of $d\nu_\theta$ is its stratified continuity in $\Theta$ can be proved by quantitative reducibility of $(\alpha,S_E^V)$. In this way, we can approximate each lost phase in  $ \Theta$ by localization phases, and prove  $d\mu_\theta^{pp}(\R)= d\nu_\theta(\R)=1$ for all phases in $ \Theta$.

Finally, we point out an interesting phenomenon based on Theorem \ref{long} and the localization result in \cite{J,bj02}: The localization phase is not sensitive. More precisely, the result in \cite{bj02} and Theorem \ref{long} imply that the localization phase of \eqref{1.2} does not sensitively depend on $V\in C^\omega(\T,\R)$ in any dimension, This phenomenon can be viewed as the robustness of localization phase, which leads us to define the robustness of Anderson localization,
\begin{Definition}
For fixed $\alpha$, $H_{V,\alpha,x}$ is said to have $C^r$ robust Anderson localization if there is a $C^r$ neighborhood $B(V)$ of $V$ and a subset $\tilde{\Theta}$, such that
$$
\bigcap_{\tilde V \in B(V)}\{ x \, | H_{\tilde{V},\alpha,x}\ \  has \ \ AL\}=\tilde\Theta,
$$
moreover $|\tilde \Theta|=1$.
\end{Definition}

\begin{Remark}
In the above definition, $\alpha$ is fixed. Similarly, one can fix $\theta$ and define another kind robustness by requiring
$$
|\bigcap_{\tilde V \in B(V)}\{ \alpha \, | H_{\tilde{V},\alpha,x}\ \  has \ \ AL\}|=1.
$$

\end{Remark}
In our forthcoming paper, we will prove that $H_{V,\alpha,x}$ with even cosine-like potential introduced in \cite{fsw,Sinai} have $C^2$ robust Anderson localization. We guess such robustness holds generally in analytic topology, however the arithmetic description of the localization phase might be more complicated. In one dimensional case, it might  relate to the acceleration defined by Avila in \cite{avila0}.

%
%
%

\section{Preliminaries}
Recall that $sl(2,\R)$ is the set of $2\times 2$ matrices  of the form
$$\left(
\begin{array}{ccc}
 x &  y+z\\
 y-z &  -x
 \end{array}\right)$$
 where $x,y,z\in \R.$ It is
 isomorphic to $su(1,1)$, the group of matrices of the form
$$\left(
\begin{array}{ccc}
 i t &  \nu\\
\bar{ \nu} &  -i t
 \end{array}\right)$$
 with $t\in \R$, $\nu\in \C$. The isomorphism between $sl(2,\R)$ and
 $su(1,1)$ is given by $B\rightarrow M B M^{-1}$ where
$$M=\frac{1}{2i}\left(
\begin{array}{ccc}
 1 &  -i\\
 1 &  i
 \end{array}\right).$$
 Direct calculation shows that
 $$M\left(
\begin{array}{ccc}
 x &  y+z\\
 y-z &  -x
 \end{array}\right)M^{-1}=\left(
\begin{array}{ccc}
 i z &  x-iy\\
x+iy &  -i z
 \end{array}\right).$$

Recall that for a bounded analytic (possibly matrix valued) function $F$ defined on $ \{ \theta |  | \Im \theta |< h \}$, let
$
\lvert F\rvert _h=  \sup_{ | \Im \theta |< h } \| F(\theta)\| $, and denote by $C^\omega_{h}(\T^d,*)$ the
set of all these $*$-valued functions ($*$ will usually denote $\R$, $sl(2,\R)$
$SL(2,\R)$).  Denote $C^{\omega}(\T^d,\R)$ by the union $\cup_{h>0}C_h^{\omega}(\T^d,\R)$.
\subsection{Quasi-periodic cocycles}

Given $A \in C^0(\T^d,{\rm SL}(2,\C))$ and rationally independent $\alpha\in\R^d$, we define the quasi-periodic \textit{cocycle} $(\alpha,A)$:
$$
(\alpha,A)\colon \left\{
\begin{array}{rcl}
\T^d \times \C^2 &\to& \T^d \times \C^2\\[1mm]
(x,v) &\mapsto& (x+\alpha,A(x)\cdot v)
\end{array}
\right.  .
$$
The iterates of $(\alpha,A)$ are of the form $(\alpha,A)^n=(n\alpha,  \mathcal{A}_n)$, where
$$
\mathcal{A}_n(x):=
\left\{\begin{array}{l l}
A(x+(n-1)\alpha) \cdots A(x+\alpha) A(x),  & n\geq 0\\[1mm]
A^{-1}(x+n\alpha) A^{-1}(x+(n+1)\alpha) \cdots A^{-1}(x-\alpha), & n <0
\end{array}\right.    .
$$
The {\it Lyapunov exponent} is defined by
$\displaystyle
L(\alpha,A):=\lim\limits_{n\to \infty} \frac{1}{n} \int_{\T^d} \ln \|\mathcal{A}_n(x)\| dx
$.

The cocycle $(\alpha,A)$ is {\it uniformly hyperbolic} if, for every $x \in \T^d$, there exists a continuous splitting $\C^2=E^s(x)\oplus E^u(x)$ such that for every $n \geq 0$,
$$
\begin{array}{rl}
|\mathcal{A}_n(x) \, v| \leq C e^{-cn}|v|, &  v \in E^s(x),\\[1mm]
|\mathcal{A}_n(x)^{-1}   v| \leq C e^{-cn}|v|, &  v \in E^u(x+n\alpha),
\end{array}
$$
for some constants $C,c>0$.
This splitting is invariant by the dynamics, i.e.,
$$A(x) E^{*}(x)=E^{*}(x+\alpha), \quad *=``s"\;\ {\rm or} \;\ ``u", \quad \forall \  x \in \T^d.$$

Assume that $A \in C^0(\T^d, {\rm SL}(2, \R))$ is homotopic to the identity. $(\alpha, A)$ induces the projective skew-product $F_A\colon \T^d \times \mathbb{S}^1 \to \T^d \times \mathbb{S}^1$ 
$$
F_A(x,w):=\left(x+\a,\, \frac{A(x) \cdot w}{|A(x) \cdot w|}\right),
$$
which is also homotopic to the identity.
Thus we can lift $F_A$ to a map $\tilde{F}_A\colon \T^d \times \R \to \T^d \times \R$ of the form $\tilde{F}_A(x,y)=(x+\alpha,y+\psi_x(y))$, where for every $x \in \T^d$, $\psi_x$ is $\Z$-periodic.
The map $\psi\colon\T^d \times \T  \to \R$ is called a {\it lift} of $A$. Let $\mu$ be any probability measure on $\T^d \times \R$ which is invariant by $\widetilde{F}_A$, and whose projection on the first coordinate is given by Lebesgue measure.
The number
$$
\rho(\alpha,A):=\int_{\T^d \times \R} \psi_x(y)\ d\mu(x,y) \ {\rm mod} \ \Z
$$
 depends  neither on the lift $\psi$ nor on the measure $\mu$, and is called the \textit{fibered rotation number} of $(\alpha,A)$ (see \cite{H,JM} for more details).

Given $\theta\in\T$, let $
R_\theta:=
\begin{pmatrix}
\cos2 \pi\theta & -\sin2\pi\theta\\
\sin2\pi\theta & \cos2\pi\theta
\end{pmatrix}$.
If $A\colon \T^d\to{\rm PSL}(2,\R)$ is homotopic to $x \mapsto R_{\frac{\langle n, x\rangle}{2}}$ for some $n\in\Z^d$,
then we call $n$ the {\it degree} of $A$ and denote it by $\deg A$.
The fibered rotation number is invariant under real conjugacies which are homotopic to the identity. More generally, if $(\alpha,A_1)$ is conjugated to $(\alpha, A_2)$, i.e., $B(x+\alpha)^{-1}A_1(x)B(x)=A_2(x)$, for some $B \colon \T^d\to{\rm PSL}(2,\R)$ with $\deg{B}=n$, then
\begin{equation}\label{rotation number}
\rho(\alpha, A_1)= \rho(\alpha, A_2)+ \frac{\langle n,\alpha \rangle}2.
\end{equation}

A cocycle $(\alpha,A)$ is called $L^2$-reducible (or $C^0$-reducible) if there exists a matrix function $B\in L^2(\T^d,SL(2,\R))$ (or $B\in C^0(\T^d,SL(2,\R))$) and a constant matrix $R$ such that
\begin{equation}\label{redu1}
B^{-1}(x+\alpha)A(x)B(x)=R,\ \  for\ \ a.e. \ \ x\in\T^d.
\end{equation}

A cocycle $(\alpha,A)$ is called $L^2$-degree 0 reducible if \eqref{redu1} holds with
$$
R=R_{\rho(\alpha,A)}.
$$

%

A typical  example is the \textit{Schr\"{o}dinger cocycles} $(\alpha,S_E^{V})$, where
$$
S_E^{V}(x):=
\begin{pmatrix}
E-V(x) & -1\\
1 & 0
\end{pmatrix},   \quad E\in\R.
$$
The Schr\"odinger cocycles are equivalent to the eigenvalue equations $H_{V, \alpha, x}u=E u$. Indeed, any formal solution $u=(u_n)_{n \in \Z}$ of $H_{V, \alpha, x}u=E u$ satisfies
$$
\begin{pmatrix}
u_{n+1}\\
u_n
\end{pmatrix}
= S_E^V(x+n\alpha) \begin{pmatrix}
u_{n}\\
u_{n-1}
\end{pmatrix},\quad \forall \  n \in \Z.
$$
The spectral properties of $H_{V,\alpha,x}$ and the dynamics of $(\alpha,S_E^V)$ are closely related by the well-known fact:
 $E\in \Sigma_{V,\alpha}$ if and only if $(\alpha,S_E^{V})$ is \textit{not} uniformly hyperbolic. Throughout the paper, we will denote $L(E)=L(\alpha,S_E^{V})$  and $\rho(E)=\rho(\alpha,S_E^{V})$ for short.

\subsection{The integrated density of states}
It is well known that the spectrum of $H_{V,\alpha,x}$ denoting by $\Sigma_{\alpha,V}$, is a compact subset of $\R$, independent of $x$ if $(1,\alpha)$ is rationally independent. The {\it integrated density of states} (IDS) $\mathcal{N}_{\alpha,V}:\R\rightarrow [0,1]$ of $H_{V,\alpha,x}$ is defined as
$$
\mathcal{N}_{\alpha,V}(E):=\int_{\T}\mu_{V,\alpha,x}(-\infty,E]dx,
$$
where $\mu_{V,\alpha,x}$ is the spectral measure of $H_{V,\alpha,x}$.

It is well known that $\rho(E)\in[0,\frac{1}{2}]$ relates to the integrated density of states $\mathcal{N}=\mathcal{N}_{\alpha,V}$ as follows:
\begin{equation}\label{relation}
\mathcal{N}(E)=1-2\rho(E).
\end{equation}

\subsection{Global theory of one-frequency Schr\"odinger operators} Let us make a short review of Avila's global theory of one frequency $SL(2,\R)$-cocycles \cite{avila0}. Suppose that $A\in C^\omega(\T, SL(2,\R))$ admits a holomorphic extension to $\{|\Im z|<h\}$. Then for $|\epsilon|<h$, we define $A_\epsilon\in C^\omega(\T,SL(2,\C))$ by $A_\epsilon(\cdot)=A(\cdot+i\epsilon)$. The cocycles which are not uniformly hyperbolic are classified into three classes: subcritical, critical, and supcritical. In particular, $(\alpha,A)$ is said to be subcritical if there exists $h>0$ such that $L(\alpha,A_\epsilon)=0$ for $|\epsilon|<h$.


A cornerstone in Avila's global theory is the ``Almost Reducibility Conjecture" (ARC), which says that $(\alpha,A)$ is almost reducible if it is subcritical. Recall that the cocycle $(\alpha,A)$ is said to be almost reducible if there exist $h_*>0$, and a sequence $B_n\in C_{h_*}^\omega(\T,PSL(2,R))$ such that $B^{-1}_n(x+\alpha)A(x)B(x)$ converges to constant uniformly in $|\Im x|<h_*$. The complete solution of ARC was recently given by Avila. 
\begin{Theorem}[Avila\cite{avila1,avila2}] \label{arc-conjecture}
Let $\alpha\in \R\backslash\Q$ and $A\in C^\omega(\T,SL(2,\R))$, $(\alpha,A)$ is almost reducible if it is subcritical.
\end{Theorem}

\subsection{Aubry duality}
Suppose that the quasi-periodic Schr\"{o}dinger operator
\begin{equation}
(H_{\lambda^{-1} V,\alpha,x}u)_n=u_{n+1}+u_{n-1}+\lambda^{-1}V(x+n\alpha)x_{n}, n\in \Z.
\end{equation}
has an analytic quasi-periodic Bloch wave $u_n=e^{2\pi i n\theta}\overline{\psi}(x+n\alpha)$ for some $\overline{\psi}\in C^\omega(\T^d,\C)$ and $\theta \in {[}0,1{)}$. It is easy to see that the Fourier coefficients of $\overline{\psi}(x)$ is an eigenfunction of  the following long range operator:
\begin{equation}
(L_{\lambda^{-1}V, \alpha,\theta}u)_m=\sum_{k\in \Z^d}V_k u_{m-k}+2\lambda\cos 2\pi(\theta+\langle m,\alpha\rangle)u_m, m\in \Z^d.
\end{equation}
$L_{\lambda^{-1}V, \alpha,\theta}$ is called the dual operator of $H_{\lambda^{-1}V,\alpha,x}$.

\section{$\mathcal{R}$-measure}
\noindent

In this section, we introduce a measure via reducibility, we call it $\mathcal{R}$-measure. We also give the relation between $\mathcal{R}$-measure and the  spectral measure.

Before introducing the $\mathcal{R}$-measure, let us first give a heuristic description on how reducibility of $(\alpha,S_E^V)$ can be used to study the pure point spectrum problem of dual model $L_{V,\alpha,\theta}$. As we introduced in Section 2.4,  reducibility of $(\alpha,S_E^V)$ can provide many Bloch waves of operator (\ref{1.1}). By Aubry dual, it thus provides many eigenfunctions for the dual operator (\ref{1.2}). The difficulty is to prove the completeness, i.e.,  those eigenfunctions form a complete basis of $\ell^2(\Z)$. Here we describe a way to prove such completeness. For a fixed $\theta$,
we define an infinite dimensional matrix $U$ whose rows are the normalized eigenfunctions provided by reducibility. If we can prove that the $\ell^2$-norm of all columns of $U$ are $1$, then $U$ defines a true unitary operator such that
$U^T L_{V,\alpha,\theta} U$ is a multiple operator, and thus (\ref{1.2}) has pure point spectrum. $\mathcal{R}$-measure is introduced based on this idea.

Assume $V\in C^0(\T^d,\R)$. We denote by
\begin{align*}
\mathcal{R}=\{E\in \Sigma_{\alpha,V}\, |\, \mbox{$(\alpha,S_E^V)$ is $C^0$-reducible and $2\rho(E)\neq \langle k,\alpha\rangle\mod\Z$, $\forall k\in\Z^d$}\}.
\end{align*}

For any $E\in \mathcal{R}$, by the definition of $C^0$-reducibility, there always exist $\bar{B}_E\in C^0(\T^d,SL(2,\R))$ and $A_E\in SL(2,\R)$ such that
\begin{equation}\label{redsch}
\bar{B}_{E}^{-1}(x+\alpha)S_{E}^V(x)\bar{B}_{E}(x)=A_{E},
\end{equation}
by \eqref{rotation number}, we have
\begin{align}\label{f1}
\rho(\alpha,A_E)=\rho(E)-\langle \ell_E,\alpha\rangle,
\end{align}
where $\ell_E=\frac{\deg{\bar{B}_E}}{2}$\footnote{We always have $\ell_E\in\Z^d$ since $\bar{B}_E\in C^0(\T^d,SL(2,\R))$.}.

Furthermore, for any $E\in \mathcal{R}$, by the definition, we have $ 2\rho(E)\neq\langle k,\alpha\rangle\mod\Z$ for any $k\in\Z^d$. By \eqref{f1}, there exists $U_E\in SL(2,\C)$ such that
  $$
  U_E^{-1}A_EU_E=\begin{pmatrix}e^{2\pi i\rho(\alpha,A_E)}&0\\0&e^{-2\pi i\rho(\alpha,A_E)}\end{pmatrix},
  $$
with $2\rho(\alpha, A_E)\neq \langle k,\alpha\rangle\mod\Z$ for any $k\in\Z^d$.

Let $B_E(x)=\bar{B}_E(x)U_E=\begin{pmatrix}b_E^{11}(x)&b_E^{12}(x)\\b_E^{21}(x)&b_E^{22}(x)\end{pmatrix}$. We define a vector-valued function $u_E:\mathcal{R}\rightarrow \ell^2(\Z)$ as the following,
  $$
u_E(n)=\frac{\hat{b}_E^{11}(n+\ell_E)}{\|b_E^{11}\|_{L^2}},
  $$
where $\hat{b}_E^{11}(n+\ell_E)=\int_{\T^d}b_E^{11}(x)e^{2\pi i\langle n+\ell_E,x\rangle}dx$.

$u_E$ may not be unique since the conjugation  $B_E$ is not unique.  However, we will show that $e^{-iarg(u_E(0))}u_E$ is unique and doesn't depend on the conjugation. Here $arg(z)=\arctan \frac{y}{x}$ if $z=x+iy$.
\begin{Lemma}\label{unique}
$e^{-iarg(u_E(0))}u_E$ is well defined for all $E\in \mathcal{R}$.
\end{Lemma}
\begin{pf}
Assume that  there exist $\tilde{B}_{E},B_{E}\in C^0(\T^d, SL(2,\C))$ with same degree\footnote{If $\ell_E\neq \tilde{\ell}_E$, let $D_E(x)=B_E(x)MR_{\langle\tilde{\ell}_E-\ell_E,x\rangle}M^{-1}=\begin{pmatrix}d_E^{11}(x)&d_E^{12}(x)\\d_E^{21}(x)&d_E^{22}(x)\end{pmatrix}$. It is obvious that $B_E$ and $D_E$ define the same $u_E$. So, one only need to prove the lemma for  $D_E$ and $\tilde B_E$. Thus, without loss of generality, we can assume that $\deg{B_E}=\deg{\tilde{B}_E}$} such that
\begin{equation}\label{id1}
\tilde{B}_{E}^{-1}(x+\alpha)S_{E}^V(x)\tilde{B}_{E}(x)=\left(
\begin{array}{ccc}
 e^{2\pi i\rho(\alpha,\tilde{A}_E)} &  0\\
0 &   e^{-2\pi i\rho(\alpha,\tilde{A}_E)}
 \end{array}\right),
\end{equation}
\begin{equation}\label{id2}
B_{E}^{-1}(x+\alpha)S_{E}^V(x)B_{E}(x)=\left(
\begin{array}{ccc}
 e^{2\pi i\rho(\alpha,A_E)} &  0\\
0 &   e^{-2\pi i\rho(\alpha,A_E)}
 \end{array}\right).
\end{equation}
It is obvious that  $\rho(\alpha,A_E)=\rho(\alpha,\tilde{A}_E)$ since $\deg{B_E}=\deg{\tilde{B}_E}$.\\

\eqref{id1} and \eqref{id2} imply that
\begin{align*}
&\ \ \ \ B_{E}(x+\alpha)\left(\begin{array}{ccc}
 e^{2\pi i\rho(\alpha,A_E)} &  0\\
0 &   e^{-2\pi i\rho(\alpha,A_E)}
 \end{array}\right)B^{-1}_{E}(x)\\
&=\tilde{B}_{E}(x+\alpha)\left(
\begin{array}{ccc}
 e^{2\pi i\rho(\alpha,A_E)} &  0\\
0 &   e^{-2\pi i\rho(\alpha,A_E)}
 \end{array}\right)\tilde{B}^{-1}_{E}(x),
\end{align*}
thus
\begin{align}\label{id3}
&\ \ \ \ T_E(x+\alpha)\left(\begin{array}{ccc}
 e^{2\pi i\rho(\alpha,A_E)} &  0\\
0 &   e^{-2\pi i\rho(\alpha,A_E)}
 \end{array}\right)\\ \nonumber
 &=\left(
\begin{array}{ccc}
 e^{2\pi i\rho(\alpha,A_E)} &  0\\
0 &   e^{-2\pi i\rho(\alpha,A_E)}
 \end{array}\right)T_E(x),
\end{align}
where $T_E(x)=\tilde{B}^{-1}_{E}B_{E}(x)=\begin{pmatrix}t_{E}^{11}(x)&t_{E}^{12}(x)\\t_{E}^{21}(x)&t_{E}^{22}(x)\end{pmatrix}$.

It follows that
\begin{equation}\label{id4}
e^{2\pi i\rho(\alpha,A_E)}t_E^{11}(x+\alpha)=e^{2\pi i\rho(\alpha,A_E)}t_E^{11}(x),
\end{equation}
\begin{equation}\label{id5}
e^{2\pi i\rho(\alpha,A_E)}t_E^{21}(x+\alpha)=e^{-2\pi i\rho(\alpha,A_E)}t_E^{21}(x),
\end{equation}
\begin{equation}\label{id6}
e^{-2\pi i\rho(\alpha,A_E)}t_E^{12}(x+\alpha)=e^{2\pi i\rho(\alpha,A_E)}t_E^{12}(x),
\end{equation}
\begin{equation}\label{id7}
e^{-2\pi i\rho(\alpha,A_E)}t_E^{22}(x+\alpha)=e^{-2\pi i\rho(\alpha,A_E)}t_E^{22}(x).
\end{equation}

By \eqref{id4}, we have
\begin{align*}
\sum\limits_{k\in\Z^d}\hat{t}_E^{11}(k)e^{2\pi i \langle k,\alpha\rangle}e^{2\pi i \langle k,x\rangle}=\sum\limits_{k\in\Z^d}\hat{t}_E^{11}(k)e^{2\pi i \langle k,x\rangle},
\end{align*}
thus
$$
(e^{2\pi i \langle k,\alpha\rangle}-1)\hat{t}_E^{11}(k)=0.
$$
Since $\alpha$ is rational independent, we have $\hat{t}_E^{11}(k)=0$ for $k\neq 0$, thus $t_E^{11}(x)=\hat{t}_E^{11}(0)$.

By \eqref{id5}, we have
\begin{align*}
&\sum\limits_{k\in\Z^d}e^{2\pi i \rho(\alpha,A_E)}\hat{t}_E^{21}(k)e^{2\pi i \langle k,\alpha\rangle}e^{2\pi i \langle k,x\rangle}\\
=&\sum\limits_{k\in\Z^d}e^{-2\pi i \rho(\alpha,A_E)}\hat{t}_E^{21}(k)e^{2\pi i \langle k,x\rangle},
\end{align*}
thus
\begin{equation*}
e^{2\pi i (2\rho(\alpha,A_E)+\langle k,\alpha\rangle)}\hat{t}_E^{21}(k)=0.
\end{equation*}
Since $2\rho(\alpha,A_E)\neq \langle k,\alpha\rangle\mod \Z$ for any $k\in\Z^d$. We have $\hat{t}_E^{21}(k)=0$ for all $k\in\Z^d$ which implies that $t_E^{21}(x)=0$.

Similarly, by \eqref{id6} and \eqref{id7}, we have $t_E^{12}(x)=0$ and $t_E^{22}(x)=\hat{t}_E^{22}(0)$. Thus
$$
T_E(x)=\begin{pmatrix}
\hat{t}^{11}_E(0)&0\\
0&\hat{t}^{22}_E(0)
\end{pmatrix},
$$
this implies that
$$
\tilde{B}^{-1}_{E}B_{E}(x)=\begin{pmatrix}
\hat{t}^{11}_E(0)&0\\
0&\hat{t}^{22}_E(0)
\end{pmatrix},
$$
thus $b^{11}_E(x)=\hat{t}^{11}_E(0)\tilde{b}^{11}_E(x)$, by the definition, we have $e^{-iarg(u_E(0))}u_E=e^{-iarg(\tilde{u}_E(0))}\tilde{u}_E$.
\end{pf}

We next define $E:\T\rightarrow \Sigma$ as the following:
$$
E(\theta)=\begin{cases}
\rho^{-1}(\theta)&\theta\in [0,\frac{1}{2}],\\
\rho^{-1}(1-\theta)& \theta\in (\frac{1}{2},1].
\end{cases}
$$
Since $\rho$ is increasing in the spectrum, $E(\theta)$ takes one value except the case $2\theta= \langle k,\alpha\rangle\mod \Z$ and the gap is open, and in this case $E(\theta)$ takes two values $\{E_-(\theta),E_+(\theta)\}$.

For fixed $\theta$, we also denote by
$$
\mathcal{N}_\theta=\{m|E_m(\theta):=E(T^m\theta)\in \mathcal{R}\},
$$
$$
\mathcal{E}_\theta=\{E_m(\theta)\}_{m\in \mathcal{N}_\theta}.
$$

By Lemma \ref{unique}, for any $E\in \mathcal{E}_\theta$, $e^{-iarg(u_E(0))}u_E$ is well defined. Thus, For any fixed $(\theta,n)\in \T\times\Z^d$, we can define the following measure,
\begin{Definition}[$\mathcal{R}$-measure]
\label{defR}
$\nu_{\theta,\delta_n}:\mathcal{B}\rightarrow \R$ is defined as:
$$
\nu_{\theta,\delta_n}(B)=\sum\limits_{m\in N_\theta^B}|u_{E_m(\theta)}(n+m)|^2.
$$
for all $B$ in the Borel $\sigma$-algebra $\mathcal{B}$ of $\R$, where $N_\theta^B=\{m| E_m(\theta)\in \mathcal{E}_\theta\cap B\}$. It is easy to see that $\nu_{\theta,\delta_n}$ is a measure, we call it $\mathcal{R}$-measure.
\end{Definition}

The following lemma motivated by \cite{jk2} gives the relationship between $\nu_{\theta,\delta_n}$ and  $\mu^{pp}_{\theta,\delta_n}$ (the pure point part of the spectral measures  $\mu_{\theta,\delta_n}$).
\begin{Lemma}\label{property}
For any $n\in\Z^d$ and $B\in\mathcal{B}$, we have
\begin{enumerate}
\item{$\nu_{\theta,\delta_n}(B)\leq \mu^{pp}_{\theta,\delta_n}(\mathcal{R}\cap B)$ for every $\theta$,}
\item{$\nu_{\theta,\delta_n}(B)=\mu^{pp}_{\theta,\delta_n}(\mathcal{R}\cap B)=\mu_{\theta,\delta_n}(\mathcal{R}\cap B)=|\mathcal{N}(\mathcal{R}\cap B)|$
for $a.e.$ $\theta$,}
\end{enumerate}
where $|\cdot|$ is the Lebesgue measure and $\mu_{\theta,\delta_n}$ is the spectral measure of $L_{V,\alpha,\theta}$ defined by
$$
\langle\delta_n,\chi_{B}(L_{V,\alpha,\theta})\delta_n\rangle=\int_{\R}\chi_{B} d\mu_{\theta,\delta_n}.
$$
\end{Lemma}
\begin{pf} We first prove (1).  For any $\theta\in\T$ and $E\in \mathcal{E}_\theta$, by the definition of $\mathcal{E}_\theta$, $(\alpha,S_E^V)$ is $C^0$-reducible, i.e., there exists $B_{E}\in C^0(\T^d, SL(2,\C))$ such that
\begin{equation}\label{red1}
B_{E}^{-1}(x+\alpha)S_{E}^V(x)B_{E}(x)=\left(
\begin{array}{ccc}
 e^{2\pi i\rho(\alpha,A_E)} &  0\\
0 &   e^{-2\pi i\rho(\alpha,A_E)}
 \end{array}\right).
\end{equation}
It follows that
\begin{equation}\label{f2}
b_E^{11}(x)=e^{2\pi i\rho(\alpha,A_E)}b_E^{21}(x+\alpha),
\end{equation}
\begin{equation}\label{f3}
(E-V(x))b_E^{11}(x)-b_E^{21}(x)=b_E^{11}(x+\alpha)e^{2\pi i\rho(\alpha,A_E)}.
\end{equation}
\eqref{f1}, \eqref{f2} and \eqref{f3} imply that
\begin{align*}
&(E-V(x))b_E^{11}(x)\\ \nonumber
=&b_E^{11}(x-\alpha)e^{-2\pi i(\rho(E)-\langle \ell_E,\alpha\rangle)}+b_E^{11}(x+\alpha)e^{2\pi i(\rho(E)-\langle \ell_E,\alpha\rangle)}.
\end{align*}

We denote by $z_E^{11}(x)=e^{-2\pi i\langle \ell_E,x\rangle}b_E^{11}(x)$, then one has
\begin{align}\label{4.1}
(E-V(x))z_E^{11}(x)=z_E^{11}(x-\alpha)e^{-2\pi i\rho(E)}+z_E^{11}(x+\alpha)e^{2\pi i\rho(E)}.
\end{align}
By the definition,
  $$
u_E(n)=\frac{\hat{b}_E^{11}(n+\ell_E)}{\|b_E^{11}\|_{L^2}}=\frac{\hat{z}_E^{11}(n)}{\|z_E^{11}\|_{L^2}},
  $$
Taking the Fourier expansion of \eqref{4.1}, we have
\begin{equation}\label{ei1}
\sum\limits_{k\in\Z^d}u_E(n-k)V_k+2\lambda\cos2\pi(\rho(E)+\langle n,\alpha\rangle)u_E(n)= Eu_E,
\end{equation}
i.e.  $\{u_E(n),n\in\Z^d\}$ is an normalized eigenfunction of the long-range operator $L_{V,\alpha,\rho(E)}$.

For any $\theta\in \T$ and $m\in\Z^d$, if $m\notin \mathcal{N}_\theta$, let $P_m(\theta)=0$. If $m\in \mathcal{N}_\theta$, let $P_m(\theta)$ be the spectral projection of $L_{V,\alpha,\theta}$ onto the eigenspace corresponding to $E_{m}(\theta)$. By the definition of $E_m(\theta)$ and the above argument, $u_{E_m(\theta)}(n)$ is an normalized eigenfunction of the long-range operator $L_{V,\alpha,T^m\theta}$, thus $T_{-m}u_{E_m(\theta)}(n)$\footnote{$T_{-m}$ is a translation defined by $T_{-m}u(n):=u(n+m)$.} is an normalized eigenfunction of the long-range operator $L_{V,\alpha,\theta}$. By the spectral theorem, we have
\begin{align*}
\mu^{pp}_{\theta,\delta_n}(\mathcal{R}\cap B)&\geq \sum\limits_{m\in\mathcal{N}_\theta^B}\langle P_{m}(\theta)\delta_n,\delta_n\rangle\\
&\geq\sum_{m\in N_{\theta}^B}|T_{-m}u_{E_m(\theta)}(n )|^2\\
&=\sum_{m\in N_{\theta}^B}|u_{E_m(\theta)}(n+m )|^2\\
&=\nu_{\theta,\delta_n}(B).
\end{align*}

Now, we prove (2). We first define a projection operator for any $\theta\in\T$ with $2\theta\neq \langle k,\alpha\rangle(\mod \Z)$ for any $k\in\Z^d$,
$$
P(\theta)=\sum_{m\in\mathcal{N}_\theta^B}P_m(\theta).
$$
Note that $2\theta\neq \langle k,\alpha\rangle$ for any $k\in\Z^d$, thus all these $E$'s in $\mathcal{E}_\theta$ are different and  all $P_m(\theta)$ are mutually orthogonal. It follows that  $P(\theta)$ is a projection. Moreover, we have
\begin{align*}
\int_{\T}\langle P(\theta)\delta_n, \delta_n\rangle d\theta&=\int_{\T}\sum\limits_{m\in\mathcal{N}_\theta^B}\langle P_m(\theta)\delta_n, \delta_n\rangle d\theta.
\end{align*}
By Fubini theorem, we have
\begin{align*}
\int_{\T}\sum\limits_{m\in\mathcal{N}_\theta^B}\langle P_m(\theta)\delta_n, \delta_n\rangle d\theta=\sum\limits_{m\in\Z^d}\int_{T^{-m}(\pm\rho(\mathcal{R}\cap B))}\langle P_m(\theta)\delta_n, \delta_n\rangle d\theta,
\end{align*}
it follows that
\begin{align*}
\sum\limits_{m\in\Z^d}\int_{T^{-m}(\pm\rho(\mathcal{R}\cap B))}\langle P_m(\theta)\delta_n, \delta_n\rangle d\theta&=\sum\limits_{m\in\Z^d}\int_{\pm\rho(\mathcal{R}\cap B)}\langle P_{m}(T^{-m}\theta)\delta_n, \delta_n\rangle d\theta.
\end{align*}
Since $T_{m}L_{V,\alpha,T^{-m}\theta}T_{-m}=L_{V,\alpha,\theta}$, we have
\begin{align*}
L_{V,\alpha,T^{-m}\theta}T_{-m}u_{E(\theta)}&=T_{-m}L_{V,\alpha,\theta}u_{E(\theta)}=E(\theta)T_{-m}u_{E(\theta)}\\
&=E_{m}(T^{-m}\theta)T_{-m}u_{E(\theta)}.
\end{align*}
It follows that $T_{-m}u_E(\theta)$ belongs to the range of $P_{m}(T^{-m}\theta)$, and for each $\delta_n\in\ell^2(\Z^d)$, we have
$$
\langle P_m(T^{-m}\theta)\delta_n,\delta_n\rangle\geq |\langle T_{-m}u_{E(\theta)},\delta_n\rangle|^2.
$$
This implies that
\begin{align*}
\sum\limits_{m\in\Z^d}\int_{\pm\rho(\mathcal{R}\cap B)}\langle P_{m}(T^{-m}\theta)\delta_n, \delta_n\rangle d\theta&\geq\sum\limits_{m\in\Z^d}\int_{\pm\rho(\mathcal{R}\cap B)}|\langle T_{-m}u_{E(\theta)}, \delta_n\rangle|^2 d\theta\\
&=\int_{\pm\rho(\mathcal{R}\cap B)}\sum\limits_{m\in\Z^d}|\langle T_mu_{E(\theta)}, \delta_n\rangle|^2 d\theta.
\end{align*}
Since $u_{E(\theta)}$ is a normalized eigenfunction, i.e., $\sum\limits_{m\in\Z^d}|\langle T_mu_{E(\theta)}, \delta_n\rangle|^2=1$.  Hence we have
\begin{align*}
\int_{\T}\langle P(\theta)\delta_n, \delta_n\rangle d\theta&\geq \int_{\pm\rho(\mathcal{R}\cap B)}\sum\limits_{m\in\Z^d}|\langle T_mu_{E(\theta)}, \delta_n\rangle|^2 d\theta\\
&=2|\rho(\mathcal{R}\cap B)|=|\mathcal{N}(\mathcal{R}\cap B)|.
\end{align*}
Together with (1), we have
\begin{align}\label{con}
|\mathcal{N}(\mathcal{R}\cap B)|&\leq \int_{\T}\nu_{\theta,\delta_n}(B)d\theta\\ \nonumber
&\leq \int_{\T}\langle P(\theta)\delta_n, \delta_n\rangle d\theta\\ \nonumber
&\leq \int_{\T}\mu^{pp}_{\theta,\delta_n}(\mathcal{R}\cap B) d\theta\\ \nonumber
&\leq \int_{\T}\mu_{\theta,\delta_n}(\mathcal{R}\cap B) d\theta\leq |\mathcal{N}(\mathcal{R}\cap B)|.
\end{align}
Thus $\nu_{\theta,\delta_n}(B)=\langle P(\theta)\delta_n, \delta_n\rangle=\mu^{pp}_{\theta,\delta_n}(\mathcal{R}\cap B)=\mu_{\theta,\delta_n}(\mathcal{R}\cap B)$ for $a.e.$ $\theta$.

Since $\mu_{\theta,\delta_n}(\mathcal{R}\cap B)=\mu_{T\theta,\delta_n}(\mathcal{R}\cap B)$, by ergodicity,  $\mu_{\theta,\delta_n}(\mathcal{R}\cap B)$ is a constant for $a.e.$ $\theta\in\T$. By \eqref{con}, one has
$$
\nu_{\theta,\delta_n}(B)=\mu^{pp}_{\theta,\delta_n}(\mathcal{R}\cap B)=\mu_{\theta,\delta_n}(\mathcal{R}\cap B)=|\mathcal{N}(\mathcal{R}\cap B)|
$$
for $a.e.$ $\theta$. This finishes  the proof.
\end{pf}
Lemma \ref{property} quickly implies the following corollary.
\begin{Corollary}\label{pro}
$\chi_{\mathcal{R}}(L_{V,\alpha,\theta})$ has pure point spectrum for $a.e.$ $\theta\in\T$. Moreover, if $|\mathcal{N}(\mathcal{R})|=1$, then $L_{V,\alpha,\theta}$ has pure point spectrum for $a.e.$ $\theta\in\T$.
\end{Corollary}

Note that Corollary \ref{pro} gives a measure version of pure point spectrum. In the following, we give a criterion for arithmetic version of pure point spectrum. Similar to the criterion of Anderson localization given in \cite{ayz} and the criterion of exponential dynamical localization in expectation given in \cite{gyz},  our criterion for arithmetic version of pure point spectrum is based on good control of $\mathcal{R}$-measure, thus based on good control of eigenfunctions.

We denote by
\begin{align}\label{def1}
\mathcal{T}_N\mathcal{E}_\theta=\{E_m(\theta)\}_{m\in \mathcal{N}_\theta,|m|\leq N},
\end{align}
\begin{align}\label{def2}
\mathcal{R}_N\mathcal{E}_\theta=\{E_m(\theta)\}_{m\in \mathcal{N}_\theta, |m|> N}.
\end{align}

The following proposition gives sufficient conditions for $\mathcal{R}$-measure being a constant on a given set.
\begin{Proposition}\label{cri1}
For any set $\mathcal{A}\subset \T$ and any $B\in\mathcal{B}$,
$$
\nu_{\theta,\delta_n}(B)=|\mathcal{N}(\mathcal{R}\cap B)|,
$$
holds for any $\theta\in\mathcal{A}$ if the following three conditions hold,
\begin{enumerate}
\item{(Uniformity Condition) For any $\epsilon>0$, there exists $N(\epsilon,\mathcal{A},B)$ such that
$$
\nu_{\theta,\delta_n}(\mathcal{R}_N\mathcal{E}_\theta\cap B)\leq \epsilon
$$
holds for $N>N(\epsilon,\mathcal{A},B)$ and any $\theta\in\mathcal{A}$.}
\item{(Continuity Condition) For any $N>0$ and any $\epsilon>0$, there exists $\delta=\delta(\epsilon,N,\mathcal{A},B)$ such that
$$
|\nu_{\theta,\delta_n}(\mathcal{T}_N\mathcal{E}_\theta\cap B)-\nu_{\theta',\delta_n}(\mathcal{T}_N\mathcal{E}_{\theta'}\cap B)|\leq \epsilon
$$
holds for any $\theta,\theta'\in\mathcal{A}$ with $|\theta-\theta'|\leq \delta$.}
\item{(Density Condition) There exists a dense subset $\mathcal{G}\subset\mathcal{A}$ such that
$$
\nu_{\theta,\delta_n}(B)=|\mathcal{N}(\mathcal{R}\cap B)|
$$
holds for any $\theta\in\mathcal{G}$.}
\end{enumerate}
\end{Proposition}
\begin{pf}
We first prove  that $\nu_{\theta,\delta_n}(B)$ is continuous on $\mathcal{A}$. For any $\epsilon>0$, by (1), there exists $N(\epsilon,\mathcal{A},B)$ such that
\begin{equation}\label{eq1}
\nu_{\theta,\delta_n}(\mathcal{R}_N\mathcal{E}_\theta\cap B)\leq \frac{\epsilon}{4}
\end{equation}
holds for $N>N(\epsilon,\mathcal{A},B)$ and any $\theta\in\mathcal{A}$.

Now we fix $N$ and $B$, by (2) there exists $\delta(\epsilon,N,A,B)$ such that
\begin{equation}\label{eq2}
|\nu_{\theta,\delta_n}(\mathcal{T}_N\mathcal{E}_\theta\cap B)-\nu_{\theta',\delta_n}(\mathcal{T}_N\mathcal{E}_{\theta'}\cap B)|\leq \frac{\epsilon}{2}
\end{equation}
holds for any $\theta,\theta'\in\mathcal{A}$ with $|\theta-\theta'|\leq \delta$.

\eqref{eq1} and \eqref{eq2} imply that for any $\epsilon>0$, there exists $\delta_0(\epsilon,\mathcal{A},B)$ such that if $|\theta-\theta'|\leq \delta_0$, we have
\begin{align*}
&\ \ \ \ |\nu_{\theta,\delta_n}(B)-\nu_{\theta',\delta_n}(B)|\\&\leq |\nu_{\theta,\delta_n}(\mathcal{T}_N\mathcal{E}_\theta\cap B)-\nu_{\theta',\delta_n}(\mathcal{T}_N\mathcal{E}_{\theta'}\cap B)|\\
&\ \ \ \ +\nu_{\theta,\delta_n}(\mathcal{R}_N\mathcal{E}_\theta\cap B)+\nu_{\theta',\delta_n}(\mathcal{R}_N\mathcal{E}_{\theta'}\cap B)\\
&\leq \frac{\epsilon}{2}+\frac{\epsilon}{2}\leq \epsilon.
\end{align*}

This proves the continuity of $\nu_{\theta,\delta_n}(B)$. Together with (3), we have $\nu_{\theta,\delta_n}(B)=|\mathcal{N}(\mathcal{R}\cap B)|$ for all $\theta\in\mathcal{A}$. Thus we finish the proof.

%
\end{pf}

We remark  that if $\mathcal{A}$ has a nice topological structure, then condition (3) in Proposition \ref{cri1} is not necessary. Recall that a closed set $\mathcal{S}\subset\R$ is called homogeneous if  there exist $\mu>0$ and $0<\sigma<diam\mathcal{S}$ such that for any $0<\epsilon<\sigma$ and any $E\in \mathcal{S}$, we have
$$
|\mathcal{S}\cap (E-\epsilon,E+\epsilon)|>\mu\epsilon.
$$
We immediately have the following corollary.
\begin{Corollary}\label{cri2}
For any homogeneous set $\mathcal{A}$ and any $B\in\mathcal{B}$,
$$
\nu_{\theta,\delta_n}(B)=|\mathcal{N}(\mathcal{R}\cap B)|,
$$
holds for all $\theta\in\mathcal{A}$ if the following two conditions hold,
\begin{enumerate}
\item{(Uniformity Condition) For any $\epsilon>0$, there exists $N(\epsilon,\mathcal{A},B)$ such that
$$
\nu_{\theta,\delta_n}(\mathcal{R}_N\mathcal{E}_\theta\cap B)\leq \epsilon
$$
holds for $N>N(\epsilon,\mathcal{A},B)$ and any $\theta\in\mathcal{A}$.}
\item{(Continuity Condition) For any $N>0$ and any $\epsilon>0$, there exists $\delta(\epsilon,N,\mathcal{A},B)$ such that
$$
|\nu_{\theta,\delta_n}(\mathcal{T}_N\mathcal{E}_\theta\cap B)-\nu_{\theta',\delta_n}(\mathcal{T}_N\mathcal{E}_{\theta'}\cap B)|\leq \epsilon
$$
holds for any $\theta,\theta'\in\mathcal{A}$ with $|\theta-\theta'|\leq \delta$.}
\end{enumerate}
\end{Corollary}
\begin{pf}
We only need to prove that (3) in Proposition \ref{cri1} holds. By conclusion (2) in Lemma \ref{property}, there exists a full measure $\mathcal{F}$ such that
$\nu_{\theta,\delta_n}(B)=\mathcal{N}(\mathcal{R}\cap B)$ for $\theta\in\mathcal{F}$. We denote by $\mathcal{G}=\mathcal{F}\cap\mathcal{A}$. For any $\theta\in\mathcal{A}$, by the definition of homogeneity, there exist $\mu>0$ and $0<\sigma<diam\mathcal{A}$ such that for any $0<\epsilon<\sigma$, we have
$$
|\mathcal{A}\cap (\theta-\epsilon,\theta+\epsilon)|>\mu\epsilon.
$$
This means $|(\theta-\epsilon,\theta+\epsilon)\cap \mathcal{G}|>0$, i.e. there exists $\theta'\in\mathcal{G}$ such that $|\theta-\theta'|\leq \epsilon$. Thus $\mathcal{G}$ is a dense subset of $\mathcal{A}$ and $\nu_{\theta,\delta_n}(B)=\mathcal{N}(\mathcal{R}\cap B)$ for $\theta\in\mathcal{G}$. By Proposition \ref{cri1},
$$
\nu_{\theta,\delta_n}(B)=|\mathcal{N}(\mathcal{R}\cap B)|,\ \ \forall \theta\in\mathcal{A}.
$$
\end{pf}

\section{Arithmetic version of Anderson localization}
In this section, we prove Theorem \ref{long} by Corollary \ref{cri2}. Denote by $\mathcal{A}_\gamma=DC_\alpha(\gamma,100\tau+d)$, we only need to prove the following theorem.
\begin{Theorem}\label{1a}
Assume that $\alpha\in DC_d(\kappa,\tau)$, $V\in C^\omega(\T^d)$ and $\mathcal{AR}=\Sigma_{\alpha,V}$. Then for any $0<\gamma<\min\{1,\frac{\kappa}{100}\}$ and any $n\in\Z^d$, we have $\nu_{\theta,\delta_n}(\mathcal{E}_\theta)=1$ for $\theta\in \mathcal{A}_\gamma$. Moreover,  all  eigenfunctions decay exponentially.
\end{Theorem}
Theorem \ref{long} is a consequence of Theorem \ref{1a}. In fact,  by the assumptions in Theorem \ref{long}, there exist $\kappa>0$ and $\tau>d-1$ such that $\alpha\in DC(\kappa,\tau)$. Moreover, since $V\in C^\omega(\T^d,\R)$ and $\mathcal{AR}=\Sigma_{\alpha,V}$, by Theorem \ref{1a}, for any $0<\gamma<\min\{1,\frac{\kappa}{100}\}$ and any $n\in\Z^d$, we have $\nu_{\theta,\delta_n}(\mathcal{E}_\theta)=1$. It follows that $\mu^{pp}_{\theta,\delta_n}(\R)=1$ for $\theta\in \mathcal{A}_\gamma$ from  (1) in Lemma \ref{property}. On the other hand, $\Theta\subset \cup_{0<\gamma<\kappa/100}\mathcal{A}_\gamma$, hence, $\mu^{pp}_{\theta,\delta_n}(\R)=1$ for $\theta\in \Theta$.  This implies that $L_{V,\alpha,\theta}$ has pure point spectrum for $\theta\in \Theta$. Together with the fact that all  eigenfunctions  decay exponentially, $L_{V,\alpha,\theta}$ has Anderson localization for $\theta\in \Theta$. \\

Now we prove Theorem \ref{1a}. It is easy to see Theorem \ref{1a} follows from  the following four conditions.
\begin{enumerate}
\item{$\mathcal{A}_\gamma$ is homogeneous for $0<\gamma<\min\{1,\frac{\kappa}{100}$\}.}
\item{$\cup_{\gamma>0}\mathcal{A}_\gamma\subset\pm\rho(\mathcal{R})$ and $u_{E(\theta)}$ decay exponentially for $\theta\in\cup_{\gamma>0}\mathcal{A}_\gamma$.}
\item{For any $\epsilon>0$, there exists $N(\epsilon,\mathcal{A}_\gamma)$ such that
$$
\nu_{\theta,\delta_n}(\mathcal{R}_N\mathcal{E}_\theta)\leq \epsilon
$$
holds for $N>N(\epsilon,\mathcal{A}_\gamma)$ and any $\theta\in\mathcal{A}_\gamma$.}
\item{For any $N>0$ and any $\epsilon>0$, there exists $\delta(\epsilon,N,\mathcal{A}_\gamma)$ such that
$$
|\nu_{\theta,\delta_n}(\mathcal{T}_N\mathcal{E}_\theta)-\nu_{\theta',\delta_n}(\mathcal{T}_N\mathcal{E}_{\theta'})|\leq \epsilon
$$
holds for any $\theta,\theta'\in\mathcal{A}_\gamma$ with $|\theta-\theta'|\leq \delta$.}
\end{enumerate}
In fact, by (1), (3), (4) and Corollary \ref{cri2}, we have $\nu_{\theta,\delta_n}(\mathcal{E}_\theta)=|\mathcal{N}(\mathcal{R})|$ for $\theta\in\mathcal{A}_\gamma$ with $\gamma<\min\{1,\frac{\kappa}{100}\}$, by (2) and \eqref{relation}, we have $|\mathcal{N}(\mathcal{R})|=1$ and all eigenfunctions $u_{E(\theta)}$ decay exponentially for $\theta\in \mathcal{A}_\gamma$.\\

Now we arrive at the final stage, i.e., the verification of  conditions (1)-(4).

\subsection{Verification of condition (1)}
We only need to prove the following Lemma.
\begin{Lemma}\label{diophantine}
Assume that $\alpha\in DC_d(\kappa,\tau)$ and $\gamma<\min\{1,\frac{\kappa}{100}\}$. Then $\mathcal{A}_\gamma$ is homogenous.
\end{Lemma}
\begin{pf}
Let $\tau'=100\tau+d$ and
$$
DC_{i,\alpha}(\gamma,\tau'):=\{\theta\in[\frac{i}{4},\frac{i+1}{4}):\|2\theta-\langle k,\alpha\rangle\|_{\R/\Z}\geq \frac{\gamma}{(|k|+1)^{\tau'}}, \forall k\in \Z\},
$$
then
$$
DC_\alpha(\gamma,\tau')=\bigcup\limits_{i=0}^{3}DC_{i,\alpha}(\gamma,\tau').
$$
We only need to prove the above result for $DC_{0,\alpha}(\gamma,\tau')$ since it is obvious that the union of two homogeneous sets is homogeneous.

Let $\Theta_k=\{\theta\in[0,\frac{1}{4}):\|2\theta-\langle k,\alpha\rangle\|_{\R/\Z}<\frac{\gamma}{(|k|+1)^{\tau'}}\}$, then $\Theta_k$ contains at most one interval. For any $\theta\in DC_{0,\alpha}(\gamma,\tau')$ and $\sigma<\min\{2^{-100\tau},\frac{\kappa^2}{16}\cdot2^{-100\tau},\frac{\gamma}{4}\}$. Set
$$
A(\theta,\sigma)=\{k:\Theta_k\cap (\theta-\sigma,\theta+\sigma)\neq \emptyset\}.
$$
Choose $k_0=k_0(\theta,\sigma)$\footnote{Note that $k_0$ maybe not unique.} such that $|k_0|=\min\limits_{k\in A(\theta,\sigma)}|k|$. Note that for any $k\in A(\theta,\sigma)$, there exist $\theta_{k_0}\in\Theta_{k_0}$, $\theta_k\in\Theta_k$ such that
\begin{equation}\label{2}
\|2\theta_{k_0}-2\theta_k\|_{\R/\Z}\leq 4\sigma,
\end{equation}
using the fact $\gamma<\frac{\kappa}{100}$, $\tau'>\tau$ and $\alpha\in DC_d(\kappa,\tau)$, we have
\begin{align}\label{3}
\|2\theta_{k_0}-2\theta_k\|_{\R/\Z}&\geq \|\langle k_0+k,\alpha\rangle\|_{\R/\Z}-\|2\theta_{k_0}-\langle k_0,\alpha\rangle\|_{\R/\Z}-\|2\theta_k-\langle k,\alpha\rangle\|_{\R/\Z}\nonumber \\
&\geq \frac{\kappa}{(|k_0+k|)^{\tau}}-\frac{\gamma}{(|k_0|+1)^{\tau'}}-\frac{\gamma}{(|k|+1)^{\tau'}}\nonumber\\
&\geq \frac{\kappa}{(4|k_0|+1)^{\tau}}.
\end{align}
\eqref{2} and \eqref{3} imply that $4\sigma\geq \frac{\kappa}{(8|k_0|+1)^{\tau}}$, thus $|k_0|\geq \frac{(\frac{\kappa}{4})^{\frac{1}{\tau}}\sigma^{-\frac{1}{\tau}}-1}{8}$.\\
Since $\sigma<\min\{2^{-100\tau},\frac{\kappa^2}{16}\cdot2^{-100\tau}\}$, we have
$$
|k_0|>\sigma^{-\frac{1}{2\tau}}.
$$
Thus
\begin{align*}
\sum\limits_{k\in A(\theta_0,\sigma)\backslash k_0}|\Theta_k\cap(\theta_0-\sigma,\theta_0+\sigma)|&\leq \sum\limits_{k\in A(\theta_0,\sigma)\backslash k_0}|\Theta_k|\\
&\leq \sum\limits_{j>\sigma^{-\frac{1}{2\tau}}}j^d\frac{\gamma}{(4j+1)^{-\tau'}}\\
&\leq \sigma.
\end{align*}
The last inequality is due to the facts $\tau'>d+100\tau$ and $\gamma<1$.

Since $\theta_0\in DC_{0,\alpha}(\gamma,\tau')$, thus $\theta_0\notin \Theta_{k_0}$ and
$$
|\Theta_{k_0}\cap(\theta_0-\sigma,\theta_0+\sigma)|\leq \sigma,
$$
which means
$$
DC_{0,\alpha}(\gamma,\tau')\cap(\theta_0-\sigma,\theta_0+\sigma)\geq \frac{1}{2}\sigma.
$$
\end{pf}
%

\subsection{Verification of condition (2)}
We first prove the following global reducibility theorem.
\begin{Theorem}\label{glored}
Assume that  $\alpha\in {\rm DC}_d(\kappa,\tau)$, $V\in C^{\omega}(\T^d,\R)$, $\gamma>0$ and $\Sigma_{\alpha,V}=\mathcal{AR}$.  Then  if $\rho(E)\in \mathcal{A}_\gamma$, there exist $\bar{B}_{E}\in C_{h_*}^{\omega}(\T^d, PSL(2,\R))$ and $A_E\in SL(2,\R)$ such that
\begin{equation}
\bar{B}_{E}(x+\alpha)^{-1}S_E^{V}(x)\bar{B}_{E}(x)=A_E,
\end{equation}
with estimate
\begin{align}\label{trans}
\|\bar{B}_{E}\|_{h_*}\leq C(\alpha,V,d,\gamma),
\end{align} where $h_*$ is a positive constant depending on $\alpha$ and $V$.
\end{Theorem}
We need the following lemma,
\begin{Lemma}[\cite{gyz}]\label{Global to local general}
Let $\alpha\in {\rm DC}_d$, $V\in C^{\omega}(\T^d,\R)$ and $\Sigma_{\alpha,V}=\mathcal{AR}$.  There exists $h_1=h_1(\alpha,V)>0$ such that for any $\eta>0$, $E\in \Sigma_{\alpha,V}$, there exists $\Phi_{E}\in C^{\omega}(\T^d, PSL(2,\R))$ with $|\Phi_{E}|_{h_1}<\Gamma(V,\alpha,\eta)$ such that
\begin{equation}
\Phi_{E}(x+\alpha)^{-1}S_E^{V}(x)\Phi_{E}(x)=R_{\phi(E)}e^{f_E(x)},
\end{equation}
with $\|f_E\|_{h_1}<\eta$.
\end{Lemma}
\begin{pf}
The proof is exactly the same as that of Lemma 5.1 in \cite{gyz}.
\end{pf}
\noindent\textbf{Proof of Theorem \ref{glored}:} Since  $h_1(\alpha,V)$ in Lemma \ref{Global to local general} is fixed, and independent of $\eta$, thus one can always take $\eta$ small enough such that
$$\eta \leq \epsilon_*(\kappa,\tau,\tau',h_1,h_1/2,d),$$
where $\epsilon_*(A_0,\kappa,\tau,\tau',h,\tilde{h},d)$ is the constant defined in  Proposition \ref{reducibility}.
Note that by Remark \ref{uniformcons}, the constant
$\epsilon_*$ given by Proposition \ref{reducibility} can be taken uniformly with respect to $R_{\phi} \in {\rm SO}(2,\R)$.
By Lemma \ref{Global to local general}, there exists $\Phi_{E}\in C_{h_1}^\omega(\T^d, PSL(2,\R))$ such that
$$\Phi_{E}^{-1}(x+\alpha)S_{E}^{V}(x)\Phi_{E}(x)=R_{\phi(E)} e^{f_{E}(x)}.$$
By footnote 5 of \cite{avila1}, $|\deg{\Phi_{E}}|\leq  C |\ln \Gamma|:=\Gamma_1$ for some constant $C=C(V,\alpha)>0$.
Since $\rho(E)\in \mathcal{A}_\gamma=DC_\alpha(\gamma,\tau')$, one has
$$
\rho(\alpha,\Phi_{E}^{-1}(\cdot+\alpha)S_{E}^{V}(\cdot)\Phi_{E}(\cdot))\in DC_\alpha(\gamma(1+\Gamma_1)^{-\tau'},\tau').
$$

By our selection, $\|f_{E}\|_{h_1}\leq \eta \leq \epsilon_*(\kappa,\tau,\tau',h_1,h_1/2,d)$. Then we can apply Proposition \ref{reducibility}, and obtain $B\in C_{h_1/2}^\omega(\T^d, PSL(2,\R))$ and $A_{E}\in SL(2,\R)$, such that
$$
B(x+\alpha)^{-1}R_{\phi(E)}e^{f_{E}(x)}B(x)=A_{E},
$$
with estimate
$$
\|B_{E}\|_{h_1/2}\leq C(\alpha,V,d,\gamma,\tau',h_1)\leq C(\alpha,V,d,\gamma).
$$
Let $\bar{B}_{E}=\Phi_{E}B$ and $h_*=\frac{h_1}{2}$, then
$$
\bar{B}_{E}^{-1}(x+\alpha)S_{E}^{V}(x)\bar{B}_{E}(x)=A_{E},
$$
with estimate
\begin{align}\label{esti1}
\|\bar{B}_{E}\|_{h_*}\leq C(\alpha,V,d,\gamma).
\end{align}\qed

We are now ready to verify condition (2). By Theorem \ref{glored}, $\cup_{\gamma>0}\mathcal{A}_\gamma\subset\pm\rho(\mathcal{R})$, thus $e^{-iarg(u_{E(\theta)}(0))}u_{E(\theta)}$ is well-defined for $\theta\in\cup_{\gamma>0}\mathcal{A}_\gamma$.  By \eqref{esti1}, $e^{-iarg(u_{E(\theta)}(0))}u_{E(\theta)}(n)$ decay exponentially since $u_{E(\theta)}(n)$ is the Fourier coefficients of $\bar{B}_{E(\theta)}U_{E(\theta)}$ for $\theta\in\cup_{\gamma>0}\mathcal{A}_\gamma$.

\subsection{Verification of condition (3)}
For any $\tilde{\gamma}>0$, $\ell\in\Z^d$, $C>0$, $0<C_{|\ell|}<1$, a normalized eigenfunction\footnote{We say $u(n)$ is normalized if $\sum_n|u(n)|^2=1$.}  $u(n)$ is said to be $(\tilde{\gamma},\ell,C,C_{|\ell|})$-good, if $$|u(n)|\leq C(e^{-\tilde{\gamma} |n|}+C_{|\ell|}e^{-\tilde{\gamma}|n+\ell|})$$ for any $n\in\Z^d$.

Since $\alpha$, $d$, $n$, $V$ are fixed,  condition (3) follows from the following theorem.
\begin{Theorem}\label{4.3}
Assume $\alpha\in DC_d(\kappa,\tau)$, $V\in C^{\omega}(\T^d,\R)$ and $\Sigma_{\alpha,V}=\mathcal{AR}$, then for any $\theta\in \mathcal{A}_\gamma$ and any $\epsilon>0$, there exists $N_0(\alpha,V,d,\gamma,n,\epsilon)$ such that for any $N>N_0$,
\begin{equation}
\nu_{\theta,\delta_n}(\mathcal{R}_N\mathcal{E}_\theta)\leq \epsilon.
\end{equation}
\end{Theorem}
The proof of Theorem \ref{4.3} is based on the following proposition proved by Ge-You-Zhou in \cite{gyz},
\begin{Proposition}[\cite{gyz}]\label{13}
Let $\alpha\in DC_d(\kappa,\tau)$, $V\in C^{\omega}(\T^d,\R)$ and $\Sigma_{\alpha,V}=\mathcal{AR}$. Then for any $\e>0$, there exist $h_1=h_1(V,\alpha)$, $C_1=C_1(V,\alpha,\e)$ with the following properties: if
$\rho(E(\theta))=\theta\in \mathcal{A}_\gamma$,  then associated with the eigenvalue $E(\theta)$,  the long range operator $L_{V,\alpha,\theta}$ has a  $$(2\pi (h_1- \frac{\e}{96}),\ell,C_1|\ln\gamma|^{4\tau'}\gamma^{-\frac{\e}{10 h_1}},\min\{1,\frac{C_1 |\ell|^{\tau'}e^{-2\pi|\ell|(h_1-\frac{\e}{96})}}{\gamma}\})$$-good eigenfunction for some $|\ell|\leq C_1|\ln\gamma|^4$.  \end{Proposition}

\noindent\textbf{Proof of Theorem \ref{4.3}:} Since $\theta\in \mathcal{A}_\gamma=DC_{\alpha}(\gamma,\tau')$,  then we have
\begin{eqnarray*}
&&\|2\theta-\langle m,\alpha\rangle-\langle k,\alpha\rangle\|_{\R/\Z} \\
&\geq& \frac{\gamma}{(|m+k|+1)^{\tau'}}\geq \frac{(1+|k|)^{-\tau'}\gamma}{(|m|+1)^{\tau'}},
\end{eqnarray*}
this implies that
\begin{align}\label{orbit}
T^k\theta\in DC_\alpha(\gamma(|k|+1)^{-\tau'},\tau').
\end{align}
By Proposition \ref{13}, for any $\e>0$, the long-range operator $L_{V,\alpha,T^k\theta}$ has a  $$(2\pi (h_1- \frac{\e}{96}),\ell,C_1|\ln\gamma'|^{4\tau'}\gamma'^{-\frac{\e}{10 h_1}},\min\{1,\frac{C_1 |\ell|^{\tau'}e^{-2\pi|\ell|(h_1-\frac{\e}{96})}}{\gamma'}\})$$-good eigenfunction for some $|\ell|\leq C_1|\ln \gamma'|^4$ where $\gamma'=\gamma(|k|+1)^{-\tau'}$.

On the other hand, by \eqref{ei1},  $\{u_{E(T^k\theta)}(n),n\in\Z^d\}$ is an normalized eigenfunction of $L_{V,\alpha,T^k\theta}$.  Let $\e=\frac{h_1}{4}$. By the definition of $(\tilde{\gamma},\ell,C_1,C_2)$-good eigenfunction, we have
\begin{align*}
|u_{E(T^k\theta)}(n+k)|&\leq C_1|\ln \gamma'|^{4\tau'}\gamma'^{-\frac{1}{40}}e^{2C_1|\ln\gamma'|^4\pi h_1}e^{-\frac{3\pi}{2}h_1|n+k|}\\
&\leq e^{-\pi h_1|k|},
\end{align*}
for $|k|\geq N_1(\alpha,V,d,\gamma,n)$, since $C_1|\ln \gamma'|^{4\tau'}\gamma'^{-\frac{1}{40}}e^{2C_1|\ln\gamma'|^4\pi h_1}$ grows sub-exponential with respect to $k$.

Thus for $N>N_1$, we have
\begin{align*}
\sum\limits_{|k|\geq N} |u_{E(T^k\theta)}(n+k)|^2\leq\sum\limits_{j\geq N} j^de^{-\pi jh_1}.
\end{align*}
Let $N_0>\max\{N_1,\frac{100|\ln\epsilon|}{\pi h_1}\}$.  By the definition of $\nu_{\theta,\delta_n}(\mathcal{R}_N\mathcal{E}_\theta)$, for $N>N_0$, we have
\begin{align*}
\nu_{\theta,\delta_n}(\mathcal{R}_N\mathcal{E}_\theta)\leq \epsilon.
\end{align*}

\subsection{Verification of condition (4)}
We only need to prove the following theorem since $\alpha$, $d$, $n$, $V$ are fixed.
\begin{Theorem}\label{continuity}
Assume $\alpha\in DC_d(\kappa,\tau)$ and $V\in C^{\omega}(\T^d,\R)$ and $\Sigma_{\alpha,V}=\mathcal{AR}$, for any $N>0$, $\epsilon>0$ and $\gamma>0$, there exists $\delta(\alpha,V,d,\gamma,n,\epsilon,N)$ such that
$$
|\nu_{\theta,\delta_n}(\mathcal{T}_N\mathcal{E}_\theta)-\nu_{\theta',\delta_n}(\mathcal{T}_N\mathcal{E}_{\theta'})|\leq \epsilon
$$
for any $\theta,\theta'\in\mathcal{A}_\gamma$ with $|\theta-\theta'|\leq \delta$.
\end{Theorem}
Theorem \ref{continuity} follows from the following two lemmas.
\begin{Lemma}\label{lem1}
For any $\epsilon>0$ and $\gamma>0$, there exists $\delta(\alpha,V,d,\gamma,\epsilon)>0$, such that if $\theta,\theta'\in \mathcal{A}_{\gamma}$ and $|\theta-\theta'|<\delta$, there exist $h_*>0$, $\bar{B}_{E(\theta)},\bar{B}_{E(\theta')}\in C_{h_*/2}^\omega(\T^d,PSL(2,\R))$ and $A_{E(\theta)},A_{E(\theta')}\in SL(2,\R)$, such that
\begin{align*}
&\ \ \ \ \bar{B}_{E(\theta)}^{-1}(x+\alpha)S_{E(\theta)}^{V}(x)\bar{B}_{E(\theta)}(x)=A_{E(\theta)},
\end{align*}
\begin{align*}
&\ \ \ \ \bar{B}_{E(\theta')}^{-1}(x+\alpha)S_{E(\theta')}^{V}(x)\bar{B}_{E(\theta')}(x)=A_{E(\theta')},
\end{align*}
with estimates
\begin{align*}
\|\bar{B}_{E(\theta)}\|_{h_*/2}\leq C(\alpha,V,d,\gamma), \ \  \forall \theta\in \mathcal{A}_{\gamma},
\end{align*}
\begin{align*}
\|\bar{B}_{E(\theta)}-\bar{B}_{E(\theta')}\|_{h_*/2}\leq \epsilon,
\end{align*}
\begin{align*}
\|A_{E(\theta)}-A_{E(\theta')}\|\leq \epsilon.
\end{align*}
\end{Lemma}

\begin{Lemma}\label{lem2}
For any $\epsilon>0$ and $\gamma>0$, there exists $\delta(\gamma,\epsilon)>0$, such that if $A,A'\in SL(2,\R)$, $\rho(A),\rho(A')\in \mathcal{A}_\gamma$ and $|A-A'|<\delta$, there exist $U,U'\in SL(2,\C)$, such that
\begin{align*}
U^{-1}AU=\begin{pmatrix}e^{2\pi i\rho(A)}&0\\ 0&e^{-2\pi i\rho(A)}\end{pmatrix},
\end{align*}
\begin{align*}
U'^{-1}A'U'=\begin{pmatrix}e^{2\pi i\rho(A')}&0\\ 0&e^{-2\pi i\rho(A')}\end{pmatrix},
\end{align*}
with estimates
\begin{align*}
\|U\|,\|U'\|\leq 2^{\tau'}\sqrt{\frac{\|A\|+\|A'\|}{\gamma}},
\end{align*}
\begin{align*}
\|U-U'\|\leq \epsilon.
\end{align*}
\end{Lemma}
We first prove Theorem \ref{continuity} assuming Lemma \ref{lem1} and Lemma \ref{lem2} hold. Then we give the proof of Lemma \ref{lem1} and Lemma \ref{lem2}.\\

\noindent
\textbf{Proof of Theorem \ref{continuity}}: For fixed $N$, since $\theta,\theta'\in\mathcal{A}_\gamma$, similar to the arguments as \eqref{orbit}, we have
$$
T^k\theta,T^k\theta'\in DC_\alpha(\gamma(1+N)^{-\tau'},\tau'),\ \  |k|\leq N.
$$
Let $\gamma_1=\gamma(1+N)^{-\tau'}$ and $\epsilon_1>0$ . By Lemma \ref{lem1}, there exists $\delta_1(\alpha.V,d,\gamma_1,\epsilon_1)$ $>0$, such that the following holds: if $T^k\theta,T^k\theta'\in \mathcal{A}_{\gamma_1}$ and $|T^k\theta-T^k\theta'|=|\theta-\theta'|<\delta_1$, then there exist $h_1>0$, $\bar{B}_{E(T^k\theta)},\bar{B}_{E(T^k\theta')}\in C_{h_1/2}^\omega(\T^d,PSL(2,\R))$ and $A_{E(T^k\theta)},A_{E(T^k\theta')}\in SL(2,\R)$, such that
\begin{align}\label{redest1}
&\ \ \ \ \bar{B}_{E(T^k\theta)}^{-1}(x+\alpha)S_{E(T^k\theta)}^{V}(x)\bar{B}_{E(T^k\theta)}(x)=A_{E(T^k\theta)},
\end{align}
\begin{align}\label{redest2}
&\ \ \ \ \bar{B}_{E(T^k\theta')}^{-1}(x+\alpha)S_{E(T^k\theta')}^{V}(x)\bar{B}_{E(T^k\theta')}(x)=A_{E(T^k\theta')},
\end{align}
with estimates
\begin{align}\label{redests1}
\|\bar{B}_{E(T^k\theta)}\|_{h_1/2},\|\bar{B}_{E(T^k\theta')}\|_{h_1/2}\leq C(\alpha,V,d,\gamma_1),
\end{align}
\begin{align}\label{redest3}
\|\bar{B}_{E(T^k\theta)}-\bar{B}_{E(T^k\theta')}\|_{h_1/2}\leq \epsilon_1,
\end{align}
\begin{align}\label{redest4}
\|A_{E(T^k\theta)}-A_{E(T^k\theta')}\|\leq \epsilon_1.
\end{align}
By \eqref{rotation number}, \eqref{redest1} and \eqref{redest2}, we have
\begin{align}\label{redest5}
\rho(\alpha, A_{E(T^k\theta)})&=T^k\theta-\frac{\langle\deg{\bar{B}_{E(T^k\theta)}},\alpha\rangle}{2},
\end{align}
\begin{align}\label{redest6}
\rho(\alpha, A_{E(T^k\theta')})&=T^k\theta'-\frac{\langle\deg{\bar{B}_{E(T^k\theta')}},\alpha\rangle}{2}.
\end{align}
by \eqref{redests1} and footnote 5 of \cite{avila1}, we have
\begin{align}\label{dest}
\deg{\bar{B}_{E(T^k\theta)}}\leq N_1(\alpha,V,d,N,\gamma).
\end{align}
Then \eqref{redest5}, \eqref{redest6} and \eqref{dest} imply
\begin{align}\label{redest7}
\rho(\alpha, A_{E(T^k\theta)})\in DC_\alpha(\gamma(1+N+N_1)^{-\tau'},\tau'),
\end{align}
\begin{align}\label{redest8}
\rho(\alpha, A_{E(T^k\theta')})\in DC_\alpha(\gamma(1+N+N_1)^{-\tau'},\tau').
\end{align}
For any $\epsilon_2>0$ and $\gamma_2=\gamma(1+N+N_1)^{-\tau'}$, by \eqref{redest7}, \eqref{redest8}, we have $\rho(A_{E(T^k\theta)}),\rho(A_{E(T^k\theta')})\in \mathcal{A}_{\gamma_2}$. By Lemma \ref{lem2}, there exists $\delta_2(\gamma_2,\epsilon_2)>0$, such that if $\|A_{E(T^k\theta)}-A_{E(T^k\theta')}\|<\delta_2$, then there exist $U_{E(T^k\theta)},U_{E(T^k\theta')}\in SL(2,\C)$, such that
\begin{align*}
&\ \ \ \ U_{E(T^k\theta)}^{-1}A_{E(T^k\theta)}U_{E(T^k\theta)}=\begin{pmatrix}e^{2\pi i\rho(\alpha,A_{E(T^k\theta)})}&0\\ 0&e^{-2\pi i\rho(\alpha,A_{E(T^k\theta)})}\end{pmatrix},
\end{align*}
\begin{align*}
&\ \ \ \ U_{E(T^k\theta')}^{-1}A_{E(T^k\theta')}U_{E(T^k\theta')}=\begin{pmatrix}e^{2\pi i\rho(\alpha,A_{E(T^k\theta')})}&0\\ 0&e^{-2\pi i\rho(\alpha,A_{E(T^k\theta')})}\end{pmatrix},
\end{align*}
with estimates
\begin{align}\label{redest9}
\|U_{E(T^k\theta)}-U_{E(T^k\theta')}\|\leq \epsilon_2,
\end{align}
\begin{align}\label{redest100}
\|U_{E(T^k\theta)}\|,\|U_{E(T^k\theta')}\|\leq 2^{\tau'}\sqrt{\frac{\|A_{E(T^k\theta)}\|+\|A_{E(T^k\theta')}\|}{\gamma}}.
\end{align}
Hence if we first choose $\epsilon_2=(\frac{\epsilon}{500(2N+1)^d})^4\cdot C^{-8}$ and then choose $\epsilon_1=\delta_2\cdot(\frac{\epsilon C^{-4}}{500(2N+1)^d})^2\cdot(2^{\tau'+1}\frac{\|A_{E(T^k\theta)}\|}{\gamma})^{-4}$, we have $\|A_{E(T^k\theta)}-A_{E(T^k\theta')}\|<\epsilon_1<\delta_2$. Then by \eqref{redests1}, \eqref{redest3}, \eqref{redest9} and \eqref{redest100}, we have
\begin{align*}
&\ \ \ \ \|B_{E(T^k\theta)}-B_{E(T^k\theta')}\|_{\frac{h_1}{2}}\\
&\leq \|\bar{B}_{E(T^k\theta)}\|_{\frac{h_1}{2}}\|U_{E(T^k\theta)}-U_{E(T^k\theta')}\|+
\|\bar{B}_{E(T^k\theta)}-\bar{B}_{E(T^k\theta')}\|_{\frac{h_1}{2}}\|U_{E(T^k\theta')}\|\\
&\leq C(\alpha,V,d,N,\gamma)\epsilon_2+\epsilon_12^{\tau'+1}\sqrt{\frac{\|A_{E(T^k\theta)}\|}{\gamma}}\\
&\leq \frac{\epsilon C^{-4}}{500(2N+1)^d}.
\end{align*}
On the one hand,
\begin{align}\label{redest10}
\big|\|b^{11}_{E(T^k\theta)}\|_{L^2}-\|b^{11}_{E(T^k\theta')}\|_{L^2}\big|&\leq \|b^{11}_{E(T^k\theta)}-b^{11}_{E(T^k\theta')}\|_{L^2}\\ \nonumber
&\leq \|B_{E(T^k\theta)}-B_{E(T^k\theta')}\|_{C^0}\\ \nonumber
&\leq \frac{\epsilon C^{-4}}{500(2N+1)^d}.
\end{align}
On the other hand, by \eqref{redest3}, we have $\ell_{E(T^k\theta)}=\ell_{E(T^k\theta')}$, thus
\begin{align}\label{redest11}
&\ \ \  \ |\hat{b}_{E(T^k\theta)}^{11}(n+\ell_{E(T^k\theta)})-\hat{b}_{E(T^k\theta')}^{11}(n+\ell_{E(T^k\theta')})|\\ \nonumber
&\leq |\int\limits_{\T}(b_{E(T^k\theta)}^{11}(x)-b_{E(T^k\theta')}^{11}(x))e^{-2\pi i(n+\ell_{E(T^k\theta)})x}dx|\\ \nonumber
&\leq\|B_{E(T^k\theta)}-B_{E(T^k\theta')}\|_{h_1/2}\\ \nonumber
&\leq \frac{\epsilon C^{-4}}{500(2N+1)^d}.
\end{align}
\eqref{redest10} and \eqref{redest11} imply for any $n\in\Z^d$,
\begin{align*}
&\ \ \ \ |u_{E(T^k\theta)}(n)-u_{E(T^k\theta')}(n)|\\
&=\big|\frac{\hat{b}_{E(T^k\theta)}^{11}(n+\ell_{E(T^k\theta)})}{\|b_{E(T^k\theta)}^{11}\|_{L^2}}-\frac{\hat{b}_{E(T^k\theta')}^{11}(n+\ell_{E(T^k\theta')})}{\|b_{E(T^k\theta')}^{11}\|_{L^2}}\big|\\
&=\frac{|\hat{b}_{E(T^k\theta)}^{11}(n+\ell_{E(T^k\theta)})\|b_{E(T^k\theta')}^{11}\|_{L^2}-\hat{b}_{E(T^k\theta')}^{11}(n+\ell_{E(T^k\theta')})\|b_{E(T^k\theta)}^{11}\|_{L^2}|}{\|b_{E(T^k\theta)}^{11}\|_{L^2}\|b_{E(T^k\theta')}^{11}\|_{L^2}}\\
&\leq \frac{\epsilon C^{-4}}{500(2N+1)^d}\frac{\|b_{E(T^k\theta)}^{11}\|_{L^2}+C}{\|b_{E(T^k\theta)}^{11}\|_{L^2}\|b_{E(T^k\theta')}^{11}\|_{L^2}}.
\end{align*}
By \eqref{f2}, we have
$$
b_{E(T^k\theta)}^{11}(x)=b_{E(T^k\theta)}^{21}(x+\alpha)e^{2\pi i\rho(\alpha,A_{E(T^k\theta)})},
$$
$$
b_{E(T^k\theta')}^{11}(x)=b_{E(T^k\theta')}^{21}(x+\alpha)e^{2\pi i\rho(\alpha,A_{E(T^k\theta')})},
$$
By the fact that  $|\det{B_{E(T^k\theta)}}|=1$, one has
\begin{align}\label{redest12}
2\|b_{E(T^k\theta)}^{11}\|_{L^2}=\|b_{E(T^k\theta)}^{11}\|_{L^2}+\|b_{E(T^k\theta)}^{21}\|_{L^2}\geq\|B_{E(T^k\theta)}\|_{C^0}^{-1},
\end{align}
\begin{align}\label{redest13}
2\|b_{E(T^k\theta')}^{11}\|_{L^2}=\|b_{E(T^k\theta')}^{11}\|_{L^2}+\|b_{E(T^k\theta')}^{21}\|_{L^2}\geq\|B_{E(T^k\theta')}\|_{C^0}^{-1}.
\end{align}
Then by \eqref{redests1}, \eqref{redest12} and \eqref{redest13}, we have
$$
\|b_{E(T^k\theta)}^{11}\|_{L^2}\|b_{E(T^k\theta')}^{11}\|_{L^2}\geq \frac{1}{4C^2}.
$$
Thus
\begin{align*}
|u_{E(T^k\theta)}(n)-u_{E(T^k\theta')}(n)|&\leq \frac{\epsilon C^{-4}}{500(2N+1)^d}\frac{\|b_{E(T^k\theta)}^{11}\|_{L^2}+C}{\|b_{E(T^k\theta)}^{11}\|_{L^2}\|b_{E(T^k\theta')}^{11}\|_{L^2}}\\
&\leq \frac{\epsilon}{100(2N+1)^d}.
\end{align*}
Let $\delta(\alpha,V,d,\gamma,\epsilon)=\min\{\delta_1,\delta_2\}$ and $|\theta-\theta'|<\delta$. By definition \ref{defR} and \eqref{def1}, one has
\begin{align*}
|\nu_{\theta,\delta_n}(\mathcal{T}_N\mathcal{E}_\theta)-\nu_{\theta',\delta_n}(\mathcal{T}_N\mathcal{E}_{\theta'})|&=|\sum\limits_{|k|\leq N}|u_{E(T^k\theta)}(n+k)|^2-\sum\limits_{|k|\leq N}|u_{E(T^k\theta')}(n+k)|^2|\\
&\leq \frac{\epsilon}{50(2N+1)^d}(2N+1)^d\leq \epsilon.
\end{align*}\qed\\

The proof of Theorem \ref{continuity} is finished. Now we give the proof of Lemma \ref{lem1} and Lemma \ref{lem2}.\\

\noindent
\textbf{Proof of Lemma \ref{lem1}:} We first consider the reducibility of cocycle $(\alpha,S_{E(\theta)}^{V})$, by Theorem \ref{glored}, there exist $\bar{B}_{E(\theta)}\in C_{h_*}^{\omega}(\T^d, PSL(2,\R))$ and $A_{E(\theta)}\in SL(2,\R)$ such that
\begin{equation}\label{pf1}
\bar{B}_{E(\theta)}(x+\alpha)^{-1}S_{E(\theta)}^{V}(x)\bar{B}_{E(\theta)}(x)=A_{E(\theta)},
\end{equation}
with estimate
$$
\|\bar{B}_{E(\theta)}\|_{h_*}\leq C(\alpha,V,d,\gamma).
$$

We now consider the reducibility of cocycle $(\alpha,S_{E(\theta')}^{V})$. Note that
\begin{align*}
&\ \ \ \ \bar{B}_{E(\theta)}^{-1}(x+\alpha)S_{E(\theta')}^{V}(x)\bar{B}_{E(\theta)}(x)\\
&=A_{E(\theta)}
+\bar{B}_{E(\theta)}^{-1}(x+\alpha)\begin{pmatrix}
E(\theta')-E(\theta)&0\\
0&0
\end{pmatrix}\bar{B}_{E(\theta)}(x).
\end{align*}
By the definition of $E(\theta)$, there exits $\delta(\epsilon')>0$ such that if $\theta,\theta'\in \mathcal{A}_\gamma$ and $|\theta-\theta'|<\delta(\epsilon')$, then $|E(\theta)-E(\theta')|<\epsilon'=C^{-100}\cdot\epsilon^8$. Thus
\begin{align}\label{estima1}
\bar{B}_{E(\theta)}^{-1}(x+\alpha)S_{E(\theta')}^{V}(x)\bar{B}_{E(\theta)}(x)&=A_{E(\theta)}+F_{E(\theta)}(x),
\end{align}
with
$$
\|F_{E(\theta)}\|_{h_*}\leq C^2|E(\theta)-E(\theta')|\leq \epsilon' C^2\leq \epsilon'^{\frac{1}{2}}.
$$
By \eqref{rotation number} and \eqref{estima1}, we have
\begin{align*}
\rho(\alpha, A_{E(\theta)}+F_{E(\theta)}(\cdot))&=\rho(\alpha, S_{E(\theta')}^{V})-\frac{\langle\deg{\bar{B}_{E(\theta)}},\alpha\rangle}{2}\\
&=\theta'-\frac{\langle\deg{\bar{B}_{E(\theta)}},\alpha\rangle}{2}.
\end{align*}
By footnote 5 of \cite{avila1}, $\deg{\bar{B}_{E(\theta)}}\leq N(\alpha,V,d,\gamma)$, thus
$$
\rho(\alpha, A_{E(\theta)}+F_{E(\theta)}(\cdot))\in DC_\alpha(\gamma(1+N)^{-\tau'},\tau').
$$
Hence there exists $\delta$ such that if $|\theta-\theta'|<\delta(\epsilon'')=\delta(\alpha,V,d,\gamma, \epsilon)$, then $|E(\theta)-E(\theta')|\leq \epsilon''=\min\{(\epsilon')^2,\frac{D_0(\gamma(1+N)^{-\tau'})^{8}}{\|A_{E(\theta)}\|^{2C_0}}(h_*/2)^{2C_0\tau'}\}$ where $C_0,D_0$ are defined in Theorem \ref{positive reducibility}, i.e.
$$
\|F_{E(\theta)}\|_{h_*}\leq C^2\epsilon''\leq (\epsilon'')^{\frac{1}{2}}\leq \frac{D_0(\gamma(1+N)^{-\tau'})^{4}}{\|A_{E(\theta)}\|^{C_0}}(h_*/2)^{C_0\tau'}.
$$
By Theorem \ref{positive reducibility}, there exists $e^{Y_{E(\theta')}}\in C_{h_*/2}^\omega(\T,SL(2,\R))$ close to the identity, such that
$$
e^{-Y_{E(\theta')}(x+\alpha)}(A_{E(\theta)}+F_{E(\theta)}(x))e^{Y_{E(\theta')}(x)}=A_{E(\theta')},
$$
with estimates $\|Y_{E(\theta')}\|_{h_*/2}\leq \epsilon''^{\frac{1}{4}}$,  $\|A_{E(\theta)}-A_{E(\theta')}\|\leq \epsilon''^{\frac{1}{4}}$.

Let $\bar{B}_{E(\theta')}=\bar{B}_{E(\theta)}e^{Y_{E(\theta')}}$. Then
\begin{equation}\label{pf2}
\bar{B}_{E(\theta')}(x+\alpha)^{-1}S_{E(\theta')}^{V}(x)\bar{B}_{E(\theta')}(x)=A_{E(\theta')},
\end{equation}
\begin{align}\label{pf3}
\|\bar{B}_{E(\theta)}-\bar{B}_{E(\theta')}\|_{\frac{h_*}{2}}\leq
\|\bar{B}_{E(\theta)}\|_{\frac{h_*}{2}}\epsilon''^{\frac{1}{4}}\leq C\epsilon'^{\frac{1}{2}}\leq \epsilon'^{\frac{1}{4}}\leq \epsilon,
\end{align}
\begin{align}\label{pf4}
\|A_{E(\theta)}-A_{E(\theta')}\|\leq \epsilon''^{\frac{1}{4}}\leq \epsilon'^{\frac{1}{2}}\leq \epsilon.
\end{align}
\eqref{pf1}, \eqref{pf2}, \eqref{pf3} and \eqref{pf4} complete the whole proof.\\

\noindent
\textbf{Proof of Lemma \ref{lem2}:} Since $\rho(A)\in \mathcal{A}_\gamma$, we have
$$
spec{A}=\{e^{2\pi i\rho(A)},e^{-2\pi i\rho(A)}\}.
$$
By Lemma 5.1 in \cite{hy}, there exists $U\in SL(2,\C)$ such that
\begin{align}\label{pf10}
U^{-1}AU=\begin{pmatrix}e^{2\pi i\rho(A)}&0\\0&e^{2\pi i\rho(A)}\end{pmatrix},
\end{align}
with estimate
$$
\|U\|\leq 2^{\tau'}\sqrt{\frac{\|A\|}{\gamma}}.
$$
Since $\|A-A'\|\leq \delta$, we can rewrite $U^{-1}A'U$ as
\begin{align*}
U^{-1}A'U&=U^{-1}(A'-A)U+\begin{pmatrix}e^{2\pi i\rho(A)}&0\\0&e^{2\pi i\rho(A)}\end{pmatrix}\\
&=\begin{pmatrix}e^{2\pi i\rho(A)}&0\\0&e^{2\pi i\rho(A)}\end{pmatrix}(I+Y),
\end{align*}
with estimate $\|Y\|\leq 2^{2\tau'}\frac{\|A\|}{\gamma}\delta$.

If $\|A'-A\|\leq \delta<\frac{1}{1000000\pi}(\frac{\gamma}{2^{\tau'}\|A\|})^{10}$, there exists $\tilde{A}\in sl(2,\C)$ with
$$
\tilde{A}=\begin{pmatrix}a&b\\c&-a\end{pmatrix},
$$
where $a=2\pi(x+iy),b,c\in\C$ and $x,y\in\R$. Such that
$$
U^{-1}A'U=exp\begin{pmatrix}a&b\\c&-a\end{pmatrix},
$$
with estimates
\begin{align}\label{et1}
|x|+|y-\rho(A)|+|b|+|c|\leq 10\delta^{\frac{1}{2}},
\end{align}
\begin{align}\label{et2}
|\rho(A')-\rho(A)|\leq \delta.
\end{align}

Since $\delta<(\frac{\gamma}{2^{\tau'}})^{10}$ and $\rho(A')\in \mathcal{A}_\gamma$, we have $|x|<10\delta^{\frac{1}{2}}<|\rho(A')|/2$, by \eqref{et1} and \eqref{et2}, we have $y\rho(A')>0$, thus
\begin{align}\label{et3}
|x+iy+i\rho(A')|\geq |y+\rho(A')|-|x|\geq \frac{1}{2}|\rho(A')|.
\end{align}

Let $P=\frac{1}{\sqrt{1+\frac{bc}{4\pi^2(x+iy+i\rho(A'))^2}}}\begin{pmatrix}1&-\frac{b}{2\pi(x+iy+i\rho(A'))}\\\frac{c}{2\pi(x+iy+i\rho(A'))}&1\end{pmatrix}$, then
$$
P^{-1}\begin{pmatrix}a&b\\c&-a\end{pmatrix}P=\begin{pmatrix}2\pi i\rho(A')&0\\0&-2\pi i\rho(A')\end{pmatrix}.
$$
By \eqref{et3}, we have
\begin{align}\label{et4}
|\frac{bc}{4\pi^2(x+iy+i\rho(A'))^2}|\leq |\frac{|b|^2+|c|^2}{\pi^2\rho(A')^2}|.
\end{align}
Since $\delta<\frac{1}{1000000\pi}(\frac{\gamma}{2^{\tau'}})^{\frac{1}{10}}$ and $\rho(A')\in \mathcal{A}_\gamma$, by \eqref{et1} we have
$$
|\frac{|b|^2+|c|^2}{\pi^2\rho(A')^2}|\leq \frac{1}{100},
$$
combining with the simple fact $|\frac{1}{\sqrt{1+x}}-1|\leq 2|x|$ for $x\leq \frac{1}{10}$, we have
$$
\|P-id\|\leq \frac{8(|b|+|c|)}{\pi\rho(A')}+|\frac{|b|^2+|c|^2}{\pi^2\rho(A')^2}|\leq \delta^{\frac{1}{4}}.
$$
For any $\epsilon>0$, let $\delta<\min\{\frac{1}{1000000\pi},(\frac{\gamma}{2^{\tau'}})^{\frac{1}{10}},\epsilon^{10}\}$ and $U'=UP$, we have
\begin{align}\label{pf11}
U'^{-1}A'U'=\begin{pmatrix}e^{2\pi i\rho(A')}&0\\0&e^{2\pi i\rho(A')}\end{pmatrix},
\end{align}
with estimate
\begin{align}\label{pf12}
\|U'-U\|\leq \|U\|\|P-id\|\leq \delta^{-\frac{1}{10}}\delta^{\frac{1}{4}}\leq \delta^{\frac{1}{10}}\leq\epsilon.
\end{align}
\eqref{pf10}, \eqref{pf11} and \eqref{pf12} complete the proof.

\section{Appendix}
We list some reducibility results for quasi-periodic $SL(2,\R)$-cocyle, one can consult \cite{ds,e,gyz} for details.
\begin{Theorem}[\cite{ds,gyz,e,hy}]\label{positive reducibility}
Let $\alpha\in DC_d(\kappa,\tau)$, $h>\tilde{h}>0$, $\tau'>d-1$, $\tau>d-1$, $\kappa>0$, $\gamma>0$,  $R\in SL(2,\R)$.   Let   $A\in C_{h}^\omega(\T^d,SL(2,\R))$ with   $\rho(\alpha,A)\in DC_\alpha(\kappa,\tau)$. Then  there exist numerical constant $C_0$, constant $D_0=D_0(\kappa,\tau,d)$,  $\epsilon=\epsilon(\tau',\tau,\kappa,\gamma,h,\tilde{h},d,R)$, such that if
$$
\|A(x)-R\|_{h}\leq \epsilon  \leq  \frac{D_0\gamma^{4}}{\|A\|^{C_0}}(h-\tilde{h})^{C_0\tau'},
$$
then there exist $B\in C_{\tilde{h}}^\omega(\T^d,SL(2,\R))$ and $\tilde{A}\in SL(2,\R)$ such that
$$
B(x+\alpha)A(\theta)B(x)^{-1}=\tilde{A},
$$
with estimates $\|B-id\|_{\tilde{h}}\leq \|A(x)-R\|_{h}^{\frac{1}{2}}$ and  $\|\tilde{A}-R\|\leq \|A(x)-R\|_{h}$ .
\end{Theorem}

\begin{Proposition}[\cite{gyz}]\label{reducibility}
For any $0<\tilde{h}<h$, $\kappa>0,\gamma>0$, $\tau>d-1$, $\tau'>d-1$.
Suppose that $\alpha\in DC_d(\kappa,\tau)$, $\rho(\alpha,A_0e^{f_0})\in DC_\alpha(\gamma,\tau')$. Then there exist $B\in C_{\tilde{h}}^\omega(\T^d, PSL(2,\R))$ and $A\in SL(2,\R)$ satisfying
$$
B^{-1}(x+\alpha)A_0e^{f_0(x)}B(x)=A,
$$  
provided that $\|f_0\|_h<\epsilon_*$ for some $\epsilon_*>0$ depending on $A_0,\kappa,\tau,\tau',h,\tilde{h},d$. In particular, $\|B\|_{\tilde{h}}\leq C(\alpha,V,d,\gamma,\tau',h,\tilde{h})$.
\end{Proposition}
\begin{Remark}\label{uniformcons}
 If $A_0$ varies in ${\rm SO}(2,\R)$, then $\epsilon_*$ can be taken uniform with respect to $A_0$.
\end{Remark}

\section*{Acknowledgement}
 J. You was partially supported by NNSF of China (11871286) and
Nankai Zhide Foundation.

\end{document}